\documentclass{article}
\usepackage[utf8]{inputenc}

\usepackage{fullpage}
\usepackage{amsmath}
\usepackage{amssymb}
\usepackage{amsfonts}
\usepackage{amsthm}
\usepackage{enumitem}
\usepackage{graphicx}
\usepackage{xcolor}
\usepackage{caption}
\usepackage{subcaption}
\usepackage[labelfont=bf]{caption}
\usepackage{comment}
\usepackage{cases}
\usepackage{xfrac}
\usepackage{float}
\usepackage{url}

\newcommand{\hcA}{\widehat{\mathcal{P}}}
\newcommand{\cA}{\mathcal{P}}
\newcommand{\bphi}{\boldsymbol{\phi}}
\newcommand{\bpsi}{\boldsymbol{\psi}}
\newcommand{\bbeta}{\boldsymbol{\beta}}
\newcommand{\bseta}{\boldsymbol{\eta}}
\newcommand{\rU}{\mathrm U}
\newcommand{\rP}{\mathrm P}
\newcommand{\rA}{\mathrm A}

\newcommand{\rY}{\mathrm Y}
\newcommand{\bu}{\mathbf u}
\newcommand{\by}{\mathbf y}
\newcommand{\blf}{\mathbf f}

\def\ww{0.45\textwidth}
\def\hh{0.35\textwidth}

\newcommand\rev[1]{#1}

\newtheorem*{remark}{Remark}

\title{Parametric Reduced Order Models for the Generalized Kuramoto--Sivashinsky Equations}

\author{
Md Rezwan Bin Mizan \thanks{Dept. of Mathematics, University of Houston, Houston, TX 77204,
mbinmiza@cougarnet.uh.edu} \quad
Maxim Olshanskii \thanks{Dept. of Mathematics, University of Houston, Houston, TX 77204, maolshanskiy@uh.edu} \quad
Ilya Timofeyev\thanks{Dept. of Mathematics, University of Houston, Houston, TX 77204, itimofey@cougarnet.uh.edu}
}
\date{}

\begin{document}

\maketitle

\begin{abstract}
The paper studies parametric Reduced Order Models (ROMs) for the Kuramoto--Sivashinsky (KS) and generalized Kuramoto--Sivashinsky (gKS) equations. We consider several POD and POD-DEIM projection ROMs with various strategies for parameter sampling and snapshot collection. The aim is to identify an approach for constructing a ROM that is efficient across a range of parameters, encompassing several regimes exhibited by the KS and gKS solutions: weakly chaotic, transitional, and quasi-periodic dynamics. We describe such an approach and demonstrate that it is essential to develop ROMs that adequately represent the short-time transient behavior of the gKS model.
\end{abstract}


\smallskip


%%%%%%%%%%%%%%%%%%%%%%%%%%%%%%%%%
\section{Introduction}
Reduced-order modeling (ROM) of systems exhibiting chaotic (turbulent) and more deterministic (laminar) behavior has been an active area of research for several decades. ROMs have been successfully applied to numerically solve various mathematical models of such systems, often formulated in terms of  partial differential equations (PDEs). A major ongoing challenge is developing ROMs that remain both effective and efficient across a range of parameter values upon which the system may depend, particularly when the solution behavior transitions from chaotic to more deterministic depending on these parameters. Solving such parametric systems is a common need in fields like uncertainty quantification and inverse modeling.

The development of ROMs for parametric fluid systems has already been addressed in the literature; see. e.g.,~\cite{karatzas2022reduced,pichi2023artificial,hess2023data,farias2023reduced,reyes2023reduced,fischer2023reduced} for the recent work on the subject. These studies primarily consider cases where the flow regime does not change under the allowed variation of parameters, remaining either laminar, transitional, or turbulent. However, there are situations, such as in shape design, where the flow's statistical behavior changes dramatically depending on the parameter values. In these cases, ROMs that perform well for laminar flows may fail to represent turbulent solutions, and vice versa.

To gain insight into the behavior of projection-based ROMs (such as POD-ROMs) under these conditions and to provide practical recommendations, this paper focuses on a prototypical model: the generalized Kuramoto–Sivashinsky (gKS) equations. The solutions to these equations exhibit various regimes, ranging from spatio-temporal chaos to quasi-periodic behavior, with the regime largely determined by the value of a critical parameter.

The development of ROMs for the KS and gKS equations has already been addressed in the literature. Existing approaches include stochastic reduced models~\cite{stinis2004stochastic,schmuck2015new,lu2017data}, POD-ROMs~\cite{lee2005reduced,zhang2019reduced,sipp2020nonlinear}, and neural network-based ROMs~\cite{menier2023cd,abadie2024topology}. However, these studies have primarily focused on recovering statistics of interest for either chaotic or deterministic solutions of the equations through reduced modeling.

The paper focuses on projection-based reduced-order models (ROMs), which construct surrogate models by projecting a high-fidelity system onto a low-dimensional, problem-dependent vector space~\cite{benner2015survey}. Prominent examples of such ROMs for dynamical systems include proper orthogonal decomposition (POD) ROMs~\cite{lumley1967structure,sirovich1987turbulence} and its variants, such as POD-DEIM~\cite{chaturantabut2010nonlinear} and balanced POD~\cite{rowley2005model}, along with PGD-ROMs~\cite{chinesta2010recent,chinesta2013proper}.
These methods derive the projection basis by leveraging information from high-fidelity solutions sampled at specific time instances and/or parameter values, commonly referred to as solution snapshots. However, constructing a general low-dimensional space that remains accurate across a wide range of parameters and over long time horizons can be highly challenging, if not infeasible. While such a space may effectively capture the training data, it often lacks predictive accuracy beyond the reference simulations.

Several approaches have been proposed to address this limitation. Partitioning strategies, introduced in~\cite{eftang2010hp,eftang2011hp,amsallem2012nonlinear}, divide the parameter domain and assign local reduced-order bases to each subdomain offline. In~\cite{amsallem2008interpolation,son2013real}, the authors suggested adapting precomputed reduced-order spaces by interpolating them for out-of-sample parameters along geodesics on the Grassmann manifold. Another approach, which utilizes the inherent tensor structure of the space-time-parameter domain to build parameter-specific local reduced spaces, was introduced in~\cite{mamonov2022interpolatory,mamonov2024tensorial}. While these techniques have shown varying degrees of success in systems with smooth parameter dependencies and relatively simple solution manifolds, it is still a challenge to extend the methodology for a parametric dynamical system exhibiting multiple regimes within the given parameter range. In this paper, we take a step in addressing this challenge and investigate the feasibility of several strategies for building projection based ROMs (POD--ROM and POD--DEIM--ROM) for the model example of the generalized Kuramoto--Sivashinsky equation.

The paper investigates several single-parameter and multi-parameter training strategies for training a reduced-order model (ROM). We examine how the selection of the training set of parameters and initial conditions affects the ROM’s ability to predict solutions across various regimes exhibited by the KS and gKS equations: weakly chaotic, transitional, and quasi-periodic dynamics. We will show that efficient prediction of solutions beyond the set of parameters and initial conditions used for training requires including additional snapshots that capture chaotic or short-time transient behaviors of the gKS model.
Additionally, the paper argues that conventional error metrics, such as the difference between FOM and ROM solutions measured in various norms, may not always be suitable for assessing ROM performance—not only in turbulent regimes but also for transitional solutions that exhibit the formation of persistent patterns.
In addition, for turbulent regimes, we examine whether ROMs reproduce the statistical behavior of solutions. To this end, we compare the power spectra computed from FOM and ROM simulations.

The remainder of the paper is organized as follows. Sections~\ref{s:Eq} and~\ref{s:ROM} review the KS and gKS equations and the parametric POD-ROM, respectively. Section~\ref{s:num}, the core of the paper, introduces our training strategies and the metrics used to assess ROM quality. Depending on the value of the problem parameter, our study covers three solution regimes: the chaotic regime (Section~\ref{s:Ch}), the transient regime with persistent pattern formation (Section~\ref{s:Tr}), and the laminar regime (Section~\ref{s:La}). Section~\ref{s:DEIM} examines the performance of the DEIM hyper-reduction technique for handling nonlinear terms in the equations. Finally, several conclusions are presented in Section~\ref{s:conc}.

%%%%%%%%%%%%%%%%%%%%%%%%%%%%%%%%%
\section{KS and gKS equations}\label{s:Eq}
Kuramoto-Sivashinsky (KS) equation 
\cite{kuramoto1976persistent,SIVASHINSKY1977}
has become a classical model for studies of spatio-temporal chaos. 
The generalized Kuramoto-Sivashinsky (gKS) equation in 1D has the following form 
\begin{equation}\label{eqn:gks}
   \frac{\partial u}{\partial t} + \frac12 \frac{(\partial u^2)}{\partial x} + \frac{\partial^2 u}{\partial x^2} + 
   \gamma \frac{\partial^3 u}{\partial x^3} +
   \frac{\partial^4 u}{\partial x^4}  = 0
\end{equation}
%\begin{equation}\label{eqn:ks}
%   \frac{\partial u}{\partial t} + \frac12 \frac{\partial u^2}{\partial x} + \frac{\partial^2 u}{\partial x^2} + \frac{\partial^4 u}{\partial x^4}  = 0
%\end{equation}
usually supplemented with periodic boundary conditions 
$u(0,t) = u(L,t)$. The KS equation is a particular case of the gKS with $\gamma=0$.
In \eqref{eqn:gks}, $L$ plays the role of the bifurcation parameter which determines the number of linearly unstable Fourier modes 
$S_{unst} = \left\lfloor \frac{L}{2\pi} \right \rfloor$, 
where $\lfloor \cdot \rfloor$
\rev{denotes the greatest integer less than or equal to that number}. The most unstable wavenumber is given by $M_{unst} = \frac{L}{2\sqrt{2}\pi}$.

The KS equation has been the subject of an extensive investigation. In particular, it was demonstrated that the KS model has a finite-dimensional attractor and, thus, is equivalent
to a finite-dimensional dynamical system \cite{constantin2012integral, foias1988inertial, nicolaenko1985some}
(see \cite{temam2012infinite} for an overview).
Moreover, it can be shown that the 
attractor is an inertial manifold \cite{jolly1990approximate, COLLET1993, temam1994}. Many other analytical and numerical studies have been performed
(e.g. \cite{kevrekidis1990back, hyman1986order, smyrlis1991predicting, munkel1996intermittency}).

The dynamics of the gKS equation also has been studied extensively 
(e.g. \cite{kawahara1983formation, manneville1985liapounov, kalliadasis2011falling,  balmforth1995solitary, chang1995interaction, ei1994equation, duprat2009liquid, tseluiko2010interaction, duprat2011wave, tseluiko2014weak, pradas2011rigorous, chang1993laminarizing}).
It has been demonstrated that parameter $\gamma$ has a strong influence on the dynamics of the gKS equation. For smaller values of $\gamma$ near zero, the system exhibits spatio-temporal chaos, similar to the KS equation. As $\gamma$ increases, the system exhibits a transition towards less chaotic or even non-chaotic dynamics
(see \cite{gotoda2015} for a detailed numerical study).
In the limit $\gamma \to \infty$ the gKS equation is equivalent to the integrable Korteweg–deVries (KdV) equation (e.g. \cite{tseluiko2010interaction}).

%%%%%%%%%%%%%%%%%%%%%%%%%%%%%%%%%
\section{Parametric POD--DEIM--ROM}\label{s:ROM}
For the \rev{spatial} discretization of \eqref{eqn:gks} we use a uniform grid on $[0,L]$ and 
apply a finite-difference discretization where nonlinear fluxes $F(u) = -u^2/2$ are discretized explicitly using $F_{i+1/2}^n = -(|u_i^n|^2 + u_i^nu_{i+1}^n + |u_{i+1}^n|^2)/6$ and linear terms are discretized implicitly. 
This discretization of the nonlinear term was used previously in \cite{zabkru65,mati00}.
The resulting system of ODEs takes the form of a parameterized dynamical systems:  
 for a given $\gamma$ from the parameter domain $\cA=[0,\Gamma]$   find $\bu = \bu(t, \gamma) : [0,T) \to \mathbb{R}^M$ solving
\begin{equation}
\label{eqn:GenericPDE}
\bu_t = \rA(\gamma)\bu+ \blf(\bu),  \quad t \in (0,T), \quad \text{and}~ \bu|_{t=0} = \bu_0.
\end{equation}
Here $\rA(\gamma)$ is an $M\times M$  matrix corresponding to the linear space derivatives in  \eqref{eqn:gks} and $\blf: \mathbb{R}^M \to \mathbb{R}^M$ stands for the discrete 
nonlinear
%'inertia' 
term; $\bu_0$ is an initial value projected on the grid.

To formulate a POD--DEIM ROM for \eqref{eqn:GenericPDE}, consider a training set of $K$ parameters sampled from the parameter 
domain, ${\hcA} := \{\widehat\gamma_1, \dots, \widehat\gamma_K \} \subset \cA$. 
Hereafter we use hats to denote parameters from the training set $\hcA$.  
At the first, offline stage of POD--DEIM, one computes through the full-order numerical simulations a collection 
of solution snapshots 
\begin{equation}
\bphi_j(\widehat\gamma_k) = \bu(t_j,\widehat\gamma_k) \in \mathbb{R}^M, 
\quad j = 1,\ldots, N, \quad k = 1,\ldots,K,
\end{equation}
and non-linear term snapshots  
\begin{equation}
\bpsi_j(\widehat\gamma_k) = \blf (\bu(t_j, \widehat\gamma_k)) 
\in \mathbb{R}^M, \quad j = 1,\ldots, N, \quad k = 1,\ldots,K,
\end{equation}
further referred to as $\bu$- and $\blf$-snapshots, respectively, at times $0\le t_1,\dots,t_N\le T$, 
and for $\widehat\gamma_k$ from the training set ${\hcA}$. For a desired reduced space dimension 
$n \rev{=n_{\rm POD}} \ll M$, one computes the reduced space basis $\{ \bu_i^{\rm pod}\}_{i=1}^{n} \subset \mathbb{R}^M$, 
referred to  as the \emph{POD basis}, such that the projection subspace 
$\mbox{span} \big\{ \bu_1^{\rm pod},\dots, \bu_n^{\rm pod}\big\}$ approximates 
the space spanned by all $\bu$-snapshots %, i.e.,  $\mbox{span}\{\bphi_j(\gamma_k)\}_{j=1,\ldots,N}^{k=1,\ldots,K}$, 
in the best possible way. 
This is achieved by assembling the matrix of all $\bu$-snapshots
\begin{equation}
\label{eqn:Phi}
\Phi_{\rm pod} = [\bphi_1(\gamma_1), \ldots, \bphi_N(\gamma_1), \ldots,
\bphi_1(\gamma_K), \ldots,\bphi_N(\gamma_K),] \in \mathbb{R}^{M \times N K}
\end{equation}
and computing its SVD. 
%\begin{equation}
%	\label{eqn:SVDa}
%\Phi_{\rm pod} = \rU \Sigma \rV^T.
%\end{equation}
Then, the POD reduced basis  $\bu_i^{\rm pod}$, $i=1,\dots,n$, consists of the first $n$ 
left singular vectors of $\Phi_{\rm pod}$. 
%, i.e., the first $n$ columns of $\rU$ or in Matlab notation $\rU_{\rm pod}=\rU_{:, 1:n}$.

Consider bow the $M\times n$ matrix $\rU_{\rm pod}=[\bu_1^{\rm pod},\ldots,\bu_n^{\rm pod}]$. At the second, online stage, the POD--ROM solution $\bu^{\rm rom}$ is found through its vector of 
coordinates $\bbeta$ in the space $\mbox{range}(\rU_{\rm pod})$, 
i.e., $\bu^{\rm rom}=\rU_{\rm pod}\bbeta$, 
which solve the projected system
\begin{equation}
\label{eqn:genericROM}
\bbeta_t = \rU_{\rm pod}^T\rA(\gamma)\rU_{\rm pod}\bbeta
+ \rU_{\rm pod}^T\blf ( \rU_{\rm pod} \bbeta),  \quad t \in (0,T), 
\mbox{~~and~~} \bbeta|_{t=0} = \rU_{\rm pod}^T\bu_0.
\end{equation}

While the projected matrix $\mathbf{U}_{\text{pod}}^T\mathbf{A}_{\gamma}\mathbf{U}_{\text{pod}}$ can be pre-computed during the offline stage, to efficiently evaluate the nonlinear term in \eqref{eqn:genericROM} during the online stage, we apply the discrete empirical interpolation method (DEIM)~\cite{chaturantabut2010nonlinear}. In this widely-used hyper-reduction technique, the nonlinear term is approximated within a lower-dimensional subspace of the space spanned by all $\blf$-snapshots.
To define the orthogonal basis $\{\mathbf{y}_i^{\text{pod}}\}_{i=1}^{n} \subset \mathbb{R}^M$ for this subspace, one takes the first $n\rev{=n_{\rm DEIM}}$ left singular vector of the matrix:
\begin{equation}
\label{eqn:Psi}
\Psi_{\rm pod} = [\bpsi_1(\gamma_1), \ldots, \bpsi_N(\gamma_1), \ldots,\bpsi_1(\gamma_K), 
\ldots,\bpsi_N(\gamma_K)] \in \mathbb{R}^{M\times NK},
 \end{equation}
consisting of all pre-computed $\blf$-snapshots.  
These $n$ first left singular vectors of 
\eqref{eqn:Psi} form the matrix 
$\rY_{\rm pod}=[ \by_1^{\rm pod},\dots, \by_{n}^{\rm pod}]$.
Then, DEIM approximates the nonlinear term of 
\eqref{eqn:GenericPDE} via
\begin{equation}
\blf (\bu)\approx \rY_{\rm pod} ( \rP^T  \rY_{\rm pod})^{-1} 
 \rP^T  \blf (\bu),
\end{equation}  
where the 'selection' matrix is defined as:
\begin{equation}
	\rP := [\mathbf{e}_{\eta_1}, \ldots, \mathbf{e}_{\eta_n}] \in \mathbb{R}^{M \times n},
\end{equation}
This matrix, $\rP$, is constructed so that for any $\blf\in\mathbb{R}^M$, the vector $\rP^T \blf$ contains $n$ entries selected from $\blf$ with indices $\bseta = [\eta_1, \ldots, \eta_n]^T \in \mathbb{R}^n$. \rev{In particular, this implies that \emph{only} those entries of $\blf$ need to be computed. Owing to the (local) definition of the discrete nonlinear fluxes, this computation involves only a small number of entries of $\bu=\rU_{\rm pod}\bbeta$ (depending on $n$), resulting in a particularly efficient treatment of the nonlinear term.}

DEIM determines $\bseta$ entirely based on the information within $\rY_{\text{pod}}$ using a greedy algorithm;  we refer to \cite{chaturantabut2010nonlinear} for further details.
 
The singular values of $\Phi_{\rm pod}$ and $\Psi_{\rm pod}$  provide information 
about the representation power of the ROMs subspaces  $\mbox{span} \big\{ \bu_1^{\rm pod},\dots, \bu_n^{\rm pod}\big\}$ 
and $\mbox{span} \big\{ \by_1^{\rm pod},\dots, \by_n^{\rm pod}\big\}$, respectively. 
In particular, the following estimate follows from the Eckart--Young--Mirsky theorem: 
\begin{equation}
\label{eqn:phipsibound}
\mbox{\small $\displaystyle  \sum_{k=1}^{K}\sum_{i=1}^{N}$ } \bigg\| \bphi_i(\gamma_k) - 
\mbox{\small $\displaystyle  \sum_{j=1}^{n}$ } \left\langle \bphi_i(\gamma_k), \bu_j^{\rm pod} \right\rangle \bu_j^{\rm pod} 
\bigg\|^2_{\ell^2} \le 
\mbox{\small $\displaystyle  \sum_{j=n+1}^{NK} $ }
\sigma_i^2 (\Phi_{\rm pod}),
\end{equation}
for representation of the solution states, and similarly for $\by^{\rm pod}_j$ and $\blf$-snapshots 
$\bpsi_i(\gamma_k)$. \rev{Here $\sigma_i (\Phi_{\rm pod})$ is the $i$th singular value of $\Phi_{\rm pod}$.}
Therefore, to determine the dimension of the reduced space, we take such minimal $n$
that 
\begin{equation}
\label{eqn:singVal}
\mbox{\small $\displaystyle  \sum_{j=n+1}^{NK} $ }
\sigma_i^2 (\Phi_{\rm pod}) \le \varepsilon \|\Phi_{\rm pod}\|_F^2
\end{equation}
with the threshold parameter $\varepsilon$. \rev{Singular values and the reduced bases $\bu_j^{\rm pod}$ and $\by_j^{\rm pod}$ are computed by performing the SVD of the snapshot matrices.}

Summarizing, the POD--DEIM ROM of \eqref{eqn:GenericPDE} takes the form 
\begin{equation} \label{POD-D-ROM}
\bbeta_t = \rU_{\rm pod}^T\rA(\gamma)\rU_{\rm pod}\bbeta + 
(\rU_{\rm pod}^T \rY_{\rm pod}) ( \rP^T  \rY_{\rm pod})^{-1}  \rP^T  
\blf ( \rU_{\rm pod}\bbeta),  \quad t \in (0,T), 
\end{equation}
with the initial condition
%\begin{equation}
$\bbeta|_{t=0} = \rU_{\rm pod}^T\bu_0$, 
%\end{equation}
for $\bbeta(t):[0,T]\to \mathbb{R}^n$ so that $\bu(t)$ is approximated by 
$\bu^{\rm rom}(t):[0,T]\to \mbox{span} \big\{ \bu_1^{\rm pod},\dots, \bu_n^{\rm pod}\big\}$, 
where $\bu^{\rm rom}(t)=\rU_{\rm pod}\bbeta(t)$. 

\begin{remark}\rm We also consider the case where the initial condition $ u_0(x,\vec{\alpha}) $ in \eqref{eqn:GenericPDE} is parametrized by a parameter vector $ \vec{\alpha} $ from a parameter domain in $ \mathbb{R}^J $. In this case, the parameter domain is sampled with $ W $ trial values for each $ \gamma_i $, where $ i = 1, \dots, K $. Therefore, the matrices $ \Phi_{\rm pod} $ and $ \Psi_{\rm pod} $ from \eqref{eqn:Phi} and \eqref{eqn:Psi} are $ M \times (NKW) $ matrices. We distinguish between the sampling of parameters $ \gamma $ and $ \vec{\alpha} $, as these parameters play very different roles in the system dynamics. For each fixed $ \vec{\alpha} $, the set of collected snapshots will be referred to as a 'trajectory'.
\end{remark}

%%%%%%%%%%%%%%%%%%%%%%%%%%%%%%%%%
\section{Numerical Results and Analysis} \label{s:num}

In this section, we present our numerical results on the performance of ROMs for the KS and gKS equations. Specifically, we focus on constructing ROMs with high predictive utility across a range of parameter values $\gamma$. This task is particularly challenging because the behavior of the gKS equation changes drastically as $\gamma$ increases (see, e.g., \cite{gotoda2015}). Notably, the chaotic regime is relatively narrow, occurring approximately within $\gamma \in [0, 0.15]$.
We demonstrate that even for larger $\gamma$, ROMs built with a sufficient number of chaotic or ``transient'' (non-chaotic regime) snapshots significantly outperform those constructed from purely non-chaotic time series. To this end, we analyze ROMs constructed from datasets specifically designed to emphasize either the chaotic (or transient) regime or the more deterministic behavior.

\begin{figure}[h]
\centering
\includegraphics[width=0.6\textwidth]{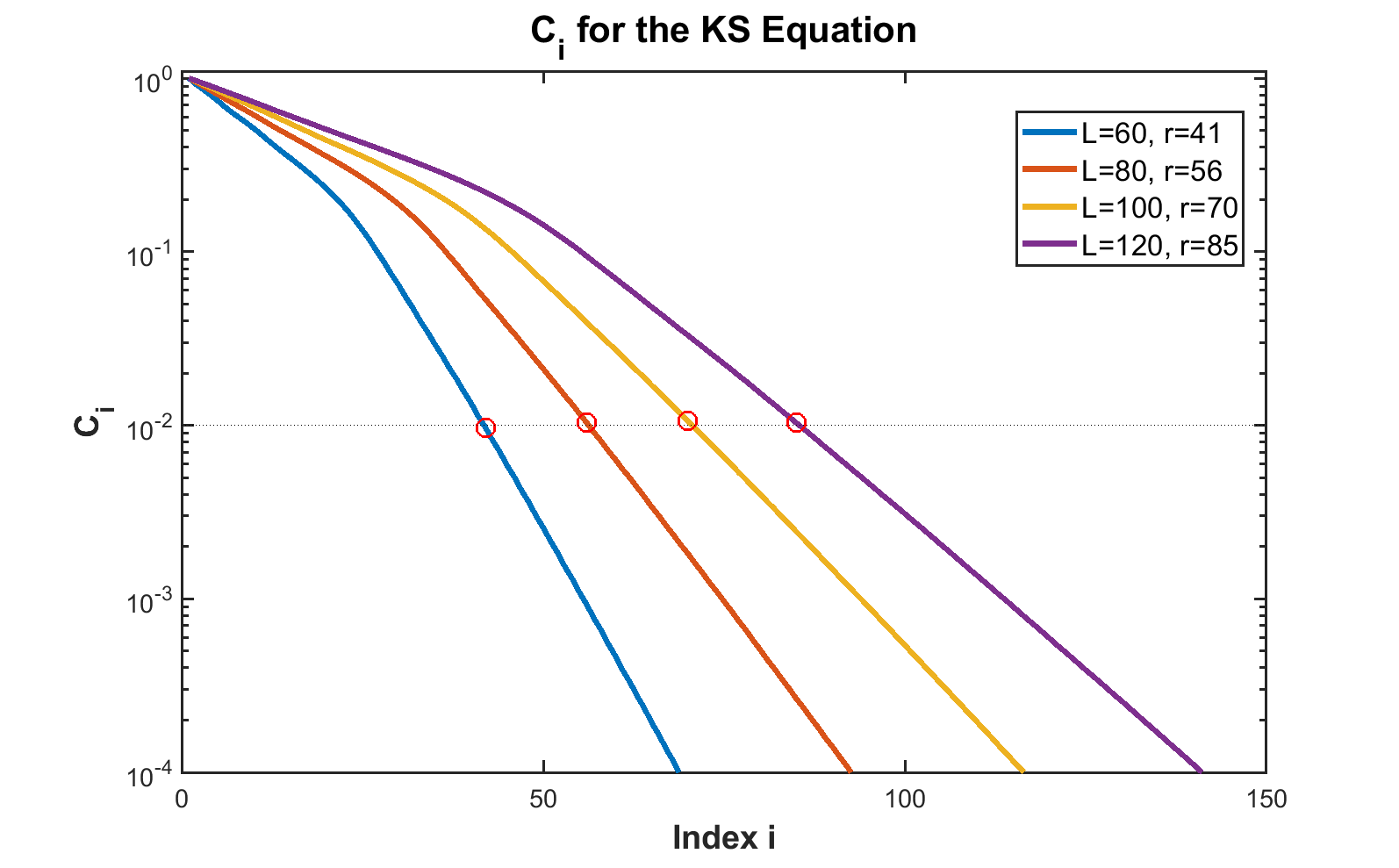}
\caption{Simulations of the KS equation. Normalized cumulative sum of singular values, $C_i$, for different values of $L$ and the corresponding optimal reduced dimension. The optimal reduced dimension is selected when $C_i$ becomes less than $10^{-2}$.}
\label{fig1}
\end{figure}

\paragraph{Parameters and Initial Conditions}
We use the following parameters in simulations of the KS and gKS equations. The number of spatial points is $M=256$, domain size $L=60$ (except for Figure \ref{fig1} and corresponding simulations).  \rev{We also performed several tests at resolution $M=512$ and compared the recovered ROM subspaces via principal angles. We found that the sine of these angles does not exceed $5\times 10^{-2}$. We therefore conclude that the spatial scales are sufficiently well resolved for the purpose of constructing the ROM.}

\rev{We use the same time step $\delta t = 0.001$ for the numerical integration of both the FOM \eqref{eqn:GenericPDE} and the ROM \eqref{POD-D-ROM}. Snapshots are not collected at every time step but at time intervals of $\Delta t = 0.5$. While several strategies for sampling snapshots are employed (see below), the total number of snapshots used for constructing the ROMs is denoted by $N^\ast$. The final integration time $T$ will vary depending on the sampling strategy and the solution dynamics. The parameter domain is $\cA=[0,10]$.}

The domain size $L=60$ corresponds to 18 unstable wavenumbers in Fourier space, and the most unstable wavenumber is 7.
For all simulations, initial conditions are generated as
\begin{equation}
u(x, 0) = \sum_{j=1}^J A_j \cos\left(\frac{2\pi j x}{L} + b_j\right),
\label{eq:initial_condition}
\end{equation}
where \(A_j \sim 0.1 \times \text{Unif}[-1, 1]\), \(b_j \sim \text{Unif}[0, 2\pi]\). 
To build our ROMs, we utilize \(J=8\), while testing the ROMs' performance and robustness, we employ a broader range of values $J=3, 8, 22$.

%%%%%%%%%%%%%%%%%%%%%%%%%%%%%%%%%
\paragraph{ROM dimensions}
% Singular Value Decomposition (SVD) is used to identify the minimum number of singular values necessary to accurately represent the data.
% %
% First, we define  \emph{cumulative sum of singular values}
% \(S_i = \sum_{j \geq i} \sigma_j\) and 
% \emph{cumulative variance ratio}
% \({C}_i= {S_i}/{S_1}\). Then we define the 
% \emph{optimal reduced dimension} as \(r=\min \{ i : {C}_i< V \}\), 
% where $V$ is a preselected threshold. In this paper, we use $V=10^{-2}$, i.e. our ROMs bases accumulate 99\% of the energy). 
%
%
\rev{As follows from \eqref{eqn:singVal}, the ROM dimensions depend on the decay rate of the singular values of the snapshot matrices. Consider the cumulative quantity \(S_i = (\sum_{j \geq i} \sigma_j^2)^{\frac12}\) which appears in \eqref{eqn:singVal}, and define the normalized quantity \({C}_i= {S_i}/{S_1}\).}
Figure~\ref{fig1} shows the decay of ${C}_i$ for four values of $L$.
To analyze the effect of the spatial domain length $L$ on reduced dimension selection,  The figure demonstrates that larger $L$ values require more singular vectors to adequately capture the system's behavior. This reveals a linear relationship between $L$ and the number of modes necessary to construct ROMs with the same threshold $\varepsilon$ as $L$ increases. Since the number of unstable wavenumbers scales linearly with $L$, the KS equation exhibits greater complexity for larger $L$, requiring a higher rank to effectively capture the dynamics of the KS and gKS models.

\subsection{Comparing Performance of Different ROMs}
From now on, we consider $L=60$.  
In our study, we construct the reduced basis for ROMs \rev{from} $N^*=20,000$ snapshots of the FOM solution. We consider several approaches for sampling these snapshots from numerical solutions. In particular, we analyze several ROMs where $\Phi_{\rm pod}$ is constructed using snapshots for a single value of $\gamma$, as well as ROMs where $\Phi_{\rm pod}$ is constructed using a combination of snapshots for multiple values of $\gamma$.

All ROMs discussed in this paper are constructed using $\varepsilon=10^{-2}$ as a threshold. The corresponding ranks vary slightly for different approaches but fall within the range $r \in [41,47]$.  

Next, we compare and contrast the following three approaches:  

\begin{itemize}
	\item \textbf{Single trajectory sampling for $\gamma=0$:}  
	Sample snapshots from a single trajectory for $\gamma=0$. The total simulation time is $T_s = 10,000$, corresponding to $N=N^*$, $K=1$, and $W=1$. The dimension of this ROM is $r=41$.  
	
	\item \textbf{Multiple trajectory sampling for $\gamma=5$:}  
	Sample snapshots from $W$ trajectories generated using different initial conditions for a single $\gamma=5$ ($K=1$). When $W$ increases, the total number of collected snapshots remains fixed at $N^*$. Consequently, trajectories become shorter, and a larger percentage of snapshots are sampled from the chaotic (or transient) regime. Specifically, we consider $W=250$, $100$, and $25$, with the total number of snapshots from a single trajectory being $N^*/W$. The resulting ranks are $46$, $47$, and $47$, respectively.  
	
	\item \textbf{Multi-parameter sampling for $\gamma$:}  
	Sample snapshots from trajectories generated using multiple values of $\gamma$. We use $W=1$ and $K=5$ different values of $\gamma$, with $N=N^*/K = 4,000$ snapshots sampled for each value. Although the number of parameters is fixed at $K=5$, we consider four different parameter domains ${\hcA}$ for sampling:  
	\begin{itemize}
		\item (i) ${\hcA}_1 = \{3, 4, 5, 7, 10\}$  
		\item (ii) ${\hcA}_2 = \{0, 4, 5, 7, 10\}$  
		\item (iii) ${\hcA}_3 = \{0, 0.3, 1, 5, 10\}$  
		\item (iv) ${\hcA}_4 = \{0, 0.2, 0.5, 0.7, 0.9\}$  
	\end{itemize}
	The corresponding ranks are $44$, $43$, $41$, and $41$, respectively.  
	
	Recall that smaller values of $\gamma$ correspond to chaotic (or transient) regimes. Therefore, ${\hcA}_1$ includes very few snapshots from the chaotic regime, while ${\hcA}_4$ includes many snapshots from the chaotic regime. We observe that the performance of the ROM improves as more chaotic snapshots are included in $\Phi_{\rm pod}$.  
\end{itemize}

The training strategies and outcome data are summarized is Table~\ref{tab:rom-strategies}.

\rev{We emphasize that in all cases the total number of snapshots $N^\ast$ is the same. Therefore, all strategies are given the same computational budget for FOM simulations. This, however, implies that the final simulation time varies. Consequently, depending on the sampling strategy, the ROMs are constructed from snapshots covering different time intervals. }

\smallskip
% >>>  ADD TABLE HERE  <
\begin{table}[ht]
\centering
\setlength{\tabcolsep}{6pt}
\renewcommand{\arraystretch}{1.25}
\footnotesize
\rev{
\begin{tabular}{|l|l|c|c|c|c|c|c|}
\hline
\textrm{Strategy} & \textrm{Training set $\widehat{\mathcal{P}}$}
& ${K}$ & ${W}$ & ${N}$ & ${T_s}$ & ${r}$ & {Figure(s)} \\
\hline

\multicolumn{8}{|l|}{\text{Single-parameter, single trajectory}} \\
\hline
ROM$_1$ & $\gamma = 0$ & 1 & 1 & 20{,}000 & 10{,}000 & 41 & 2, 4, 5, 12--14 \\
\hline

\multicolumn{8}{|l|}{\text{Single-parameter, multiple trajectories ($\gamma = 5$)}} \\
\hline
 & $\gamma = 5,\; W = 25$  & 1 & 25  & 800 & 400 & 47 & 2 \\
\cline{2-8}
ROM$_2$ & $\gamma = 5,\; W = 100$ & 1 & 100 & 200 & 100 & 47 & 2 \\
\cline{2-8}
 & $\gamma = 5,\; W = 250$ & 1 & 250 & 80 & 40 & 46 & 2, 4, 5, 7--11 \\
\hline

\multicolumn{8}{|l|}{\text{Multi-parameter ($K = 5,\; W = 1$)}} \\
\hline
 & $\widehat{\mathcal{P}}_1 = \{3,4,5,7,10\}$ & 5 & 1 & 4{,}000 & 2{,}000 & 44 & 3 \\
\cline{2-8}
 & $\widehat{\mathcal{P}}_2 = \{0,4,5,7,10\}$ & 5 & 1 & 4{,}000 & 2{,}000 & 43 & 3 \\
\cline{2-8}
 & $\widehat{\mathcal{P}}_3 = \{0,0.3,1,5,10\}$ & 5 & 1 & 4{,}000 & 2{,}000 & 41 & 3 \\
\cline{2-8}
ROM$_3$ & $\widehat{\mathcal{P}}_4 = \{0,0.2,0.5,0.7,0.9\}$ & 5 & 1 & 4{,}000 & 2{,}000 & 41 & 3, 4, 5, 7--11 \\
\hline
\end{tabular}
}
\caption{Overview of the ROM training strategies. In all 
cases, the total number of snapshots is $N^*=20{,}000$ and the snapshot time-step 
is $\Delta t = 0.5$. The number of snapshots per trajectory is $N = N^*/(KW)$, and 
the corresponding simulation time per trajectory is $T_s = N\,\Delta t = (N^*/KW)\Delta t$.
}
\label{tab:rom-strategies}
\end{table}

\smallskip

\emph{Performance assessment:} To gauge the performance of ROMs, we will focus on several statistics such as (i) Prediction time (measures the short-term prediction capability of the ROM); (ii) Power spectra (measures the statistical properties of chaotic solutions); (iii) Persistent patterns (assesses the prediction of long-term solution behavior).

\subsubsection{Averaged Prediction Time.}
First, we evaluate the performance of different reduced order models 
by comparing the prediction accuracy for several randomly chosen initial conditions, \rev{not from the training set of initial conditions}. For this purpose, we define the prediction time of the ROM for a given initial condition as 
\begin{equation}\label{eq:predT}
T_{\rm rom}(u_0) = \mathrm{arg}\max_{t\in[0,T]}\Big\{t\,:\, \sup_{t'\in[0,t]}  \frac{\|u^{\rm rom}(t')-u^{\rm fom}(t')\|_{L^2(0,L)}}{\|u^{\rm fom}(t')\|_{L^2(0,L)}}\le 10^{-1} \Big\},
\end{equation}
where $u^{\rm fom}$ solves \eqref{eqn:GenericPDE}.

\begin{figure}
\centerline{\includegraphics[width=\ww, height=\hh]{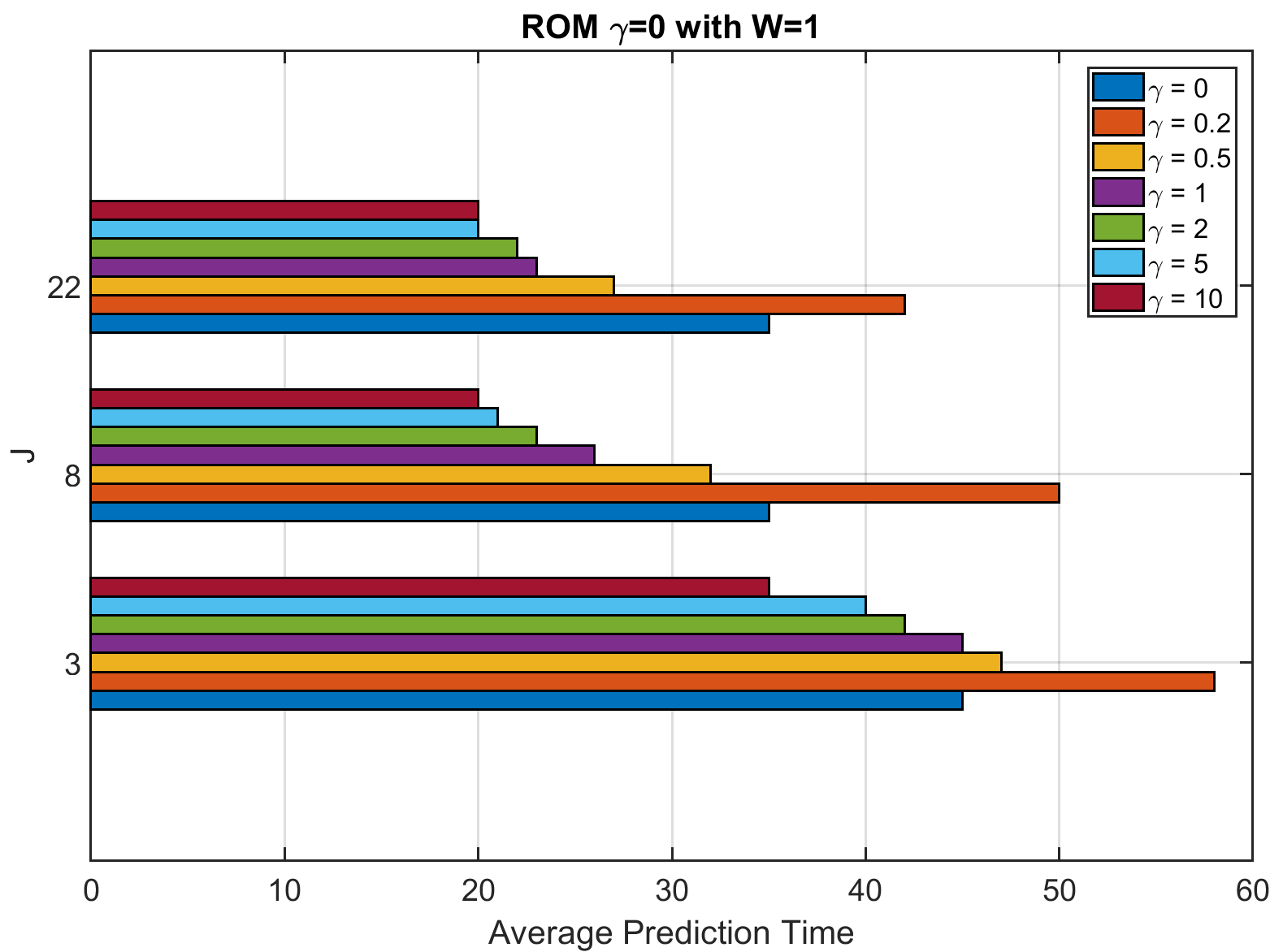}
\includegraphics[width=\ww, height=\hh]{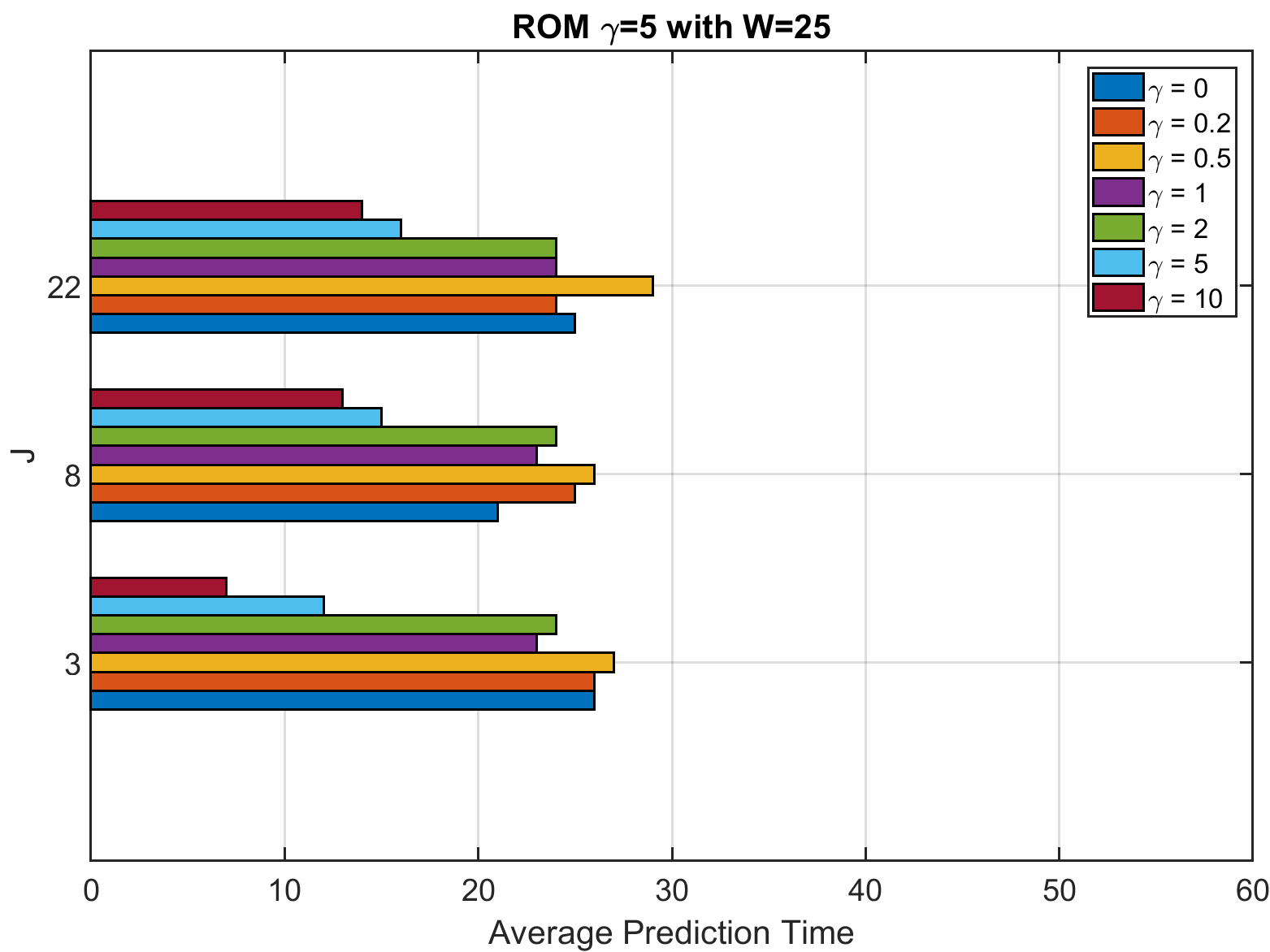}}
\vspace{0.5cm}
\centerline{\includegraphics[width=\ww, height=\hh]{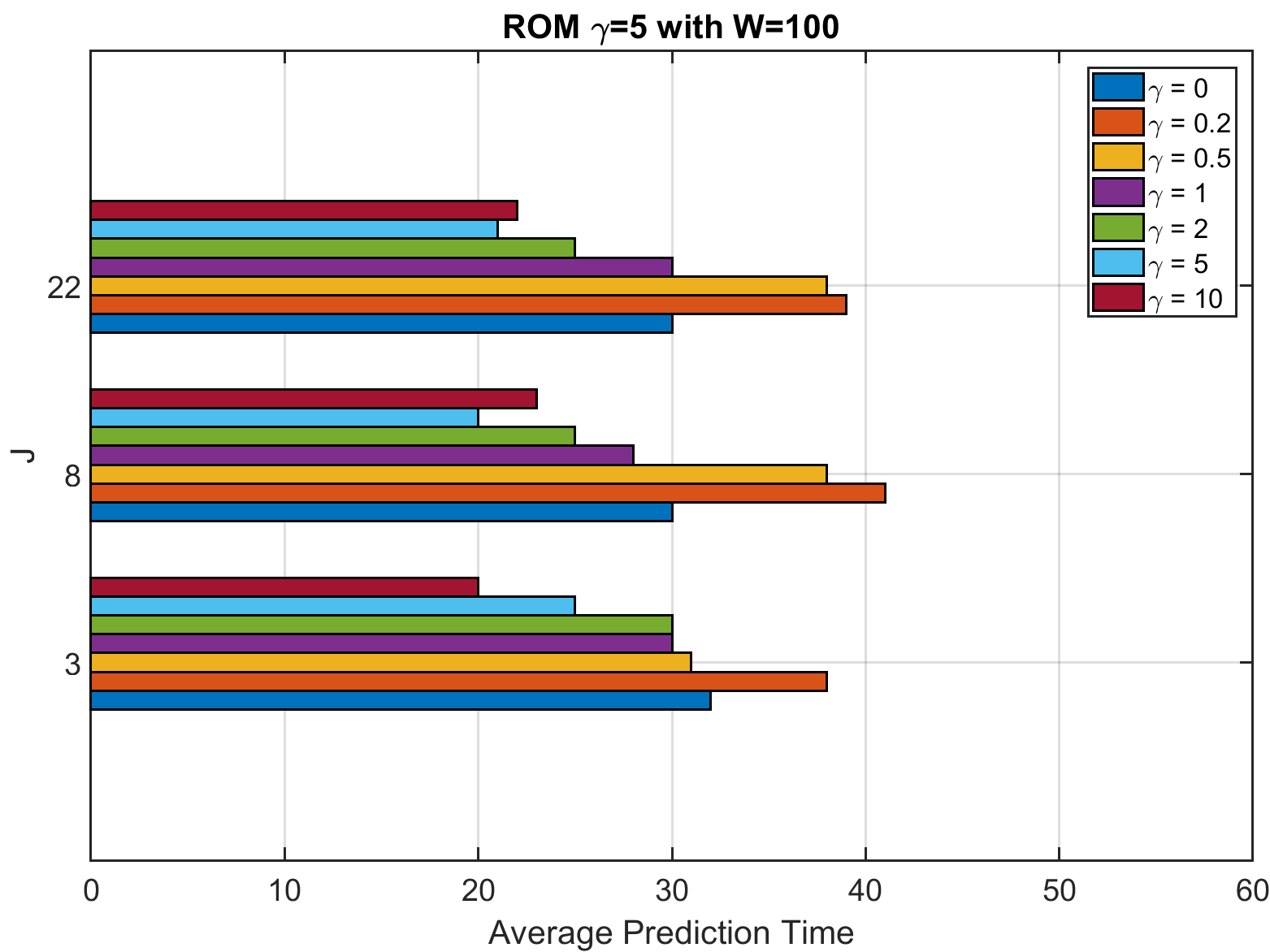}
\includegraphics[width=\ww, height=\hh]{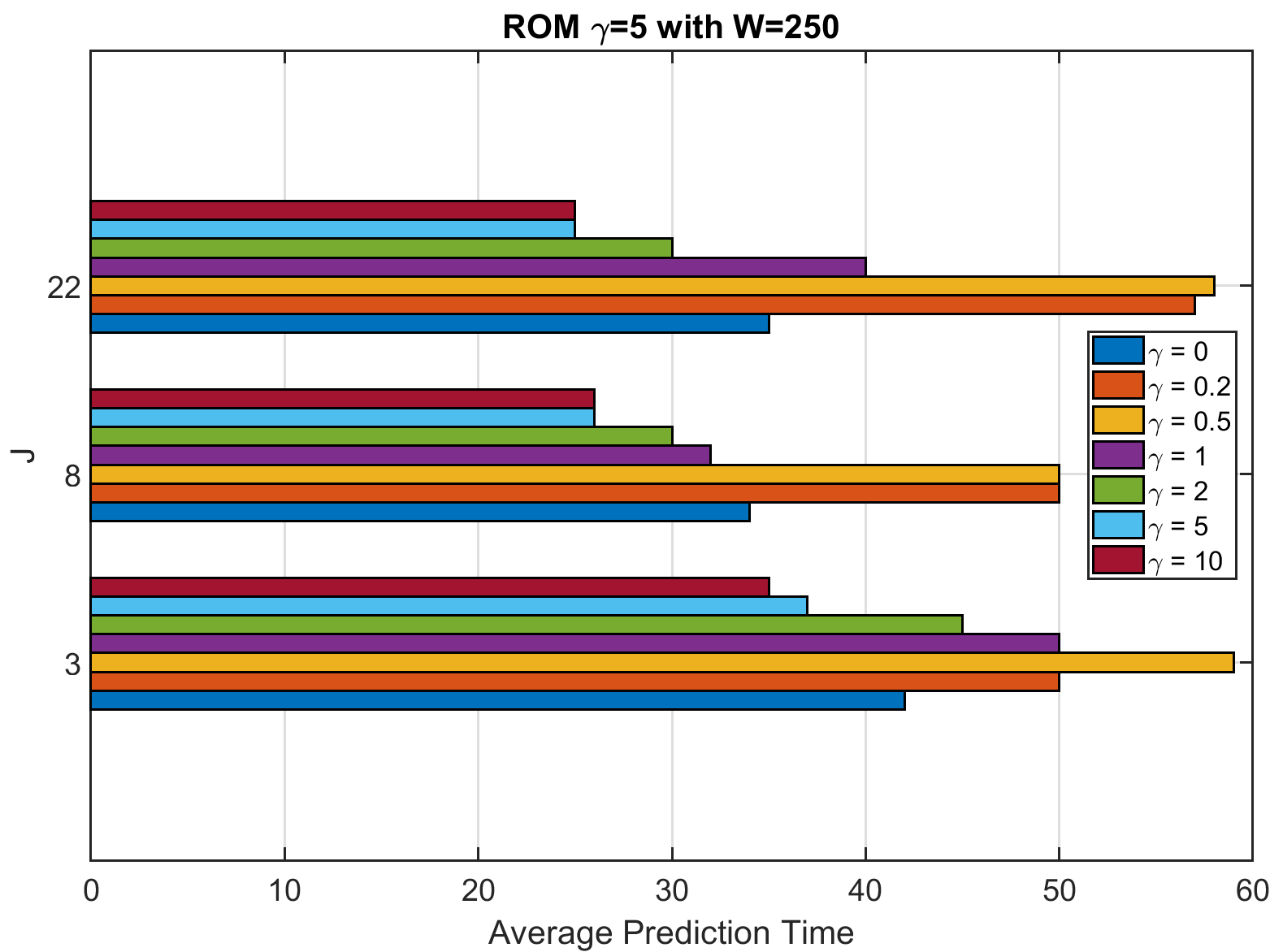}}
\caption{Averaged prediction times for trajectories of the gKS equation with different values of $\gamma$ using ROMs constructed using gKS data with a single value of $\gamma$. Upper left - ROM is constructed from simulations with $\gamma=0$ and $W=1$. Upper right, bottom left, bottom right - ROM is constructed from simulations with $\gamma=5$ with $W=25,100,250$, respectively. We consider 3 sets of different initial conditions with $J=3,8,22$, where $J$ is the number of non-zero Fourier wavenumbers at time $t=0$ (see eq. \eqref{eq:initial_condition}).}
\label{fig2}
\end{figure} 

\begin{figure}
\centerline{\includegraphics[width=\ww, height=\hh]{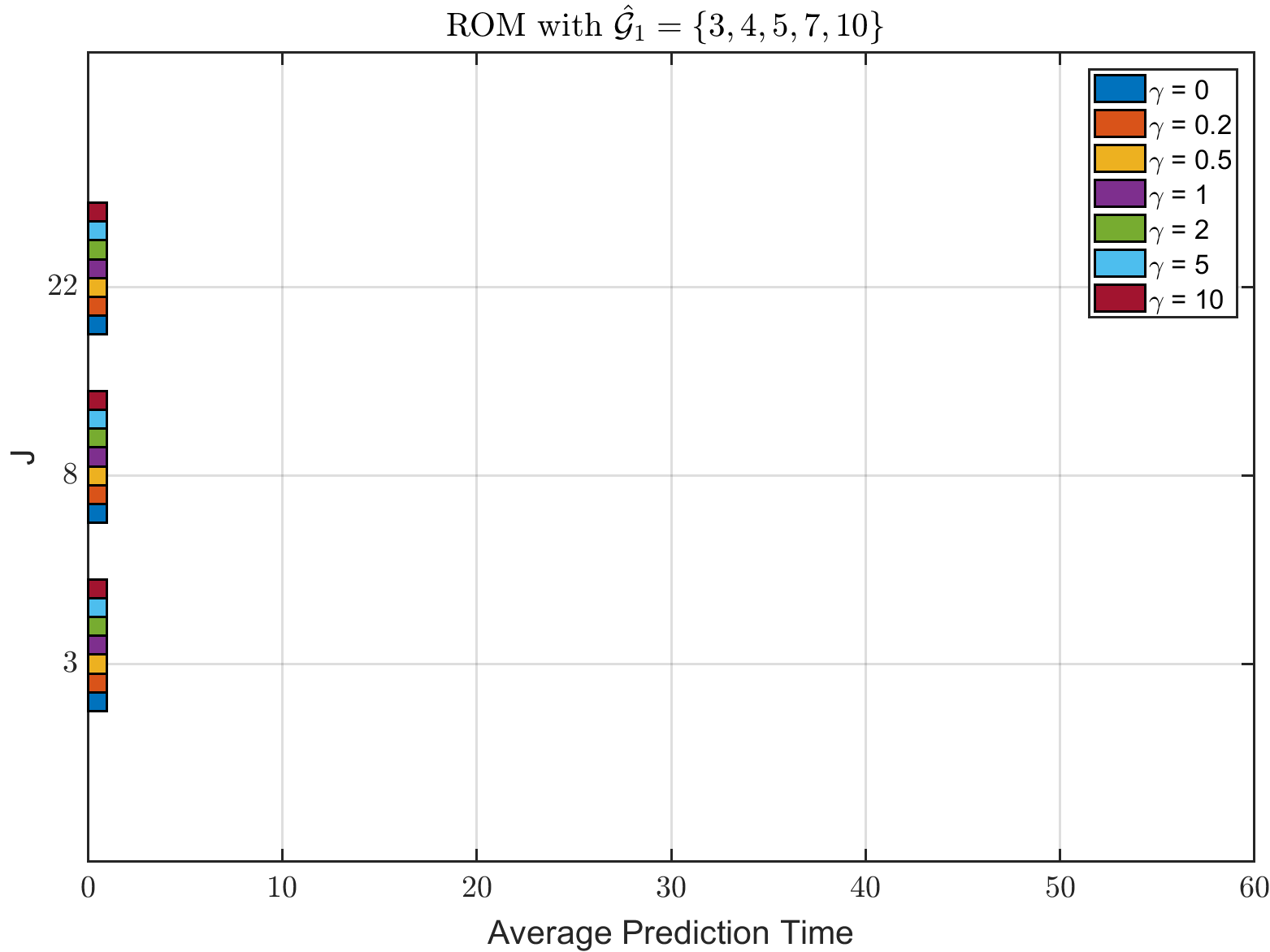}
\includegraphics[width=\ww, height=\hh]{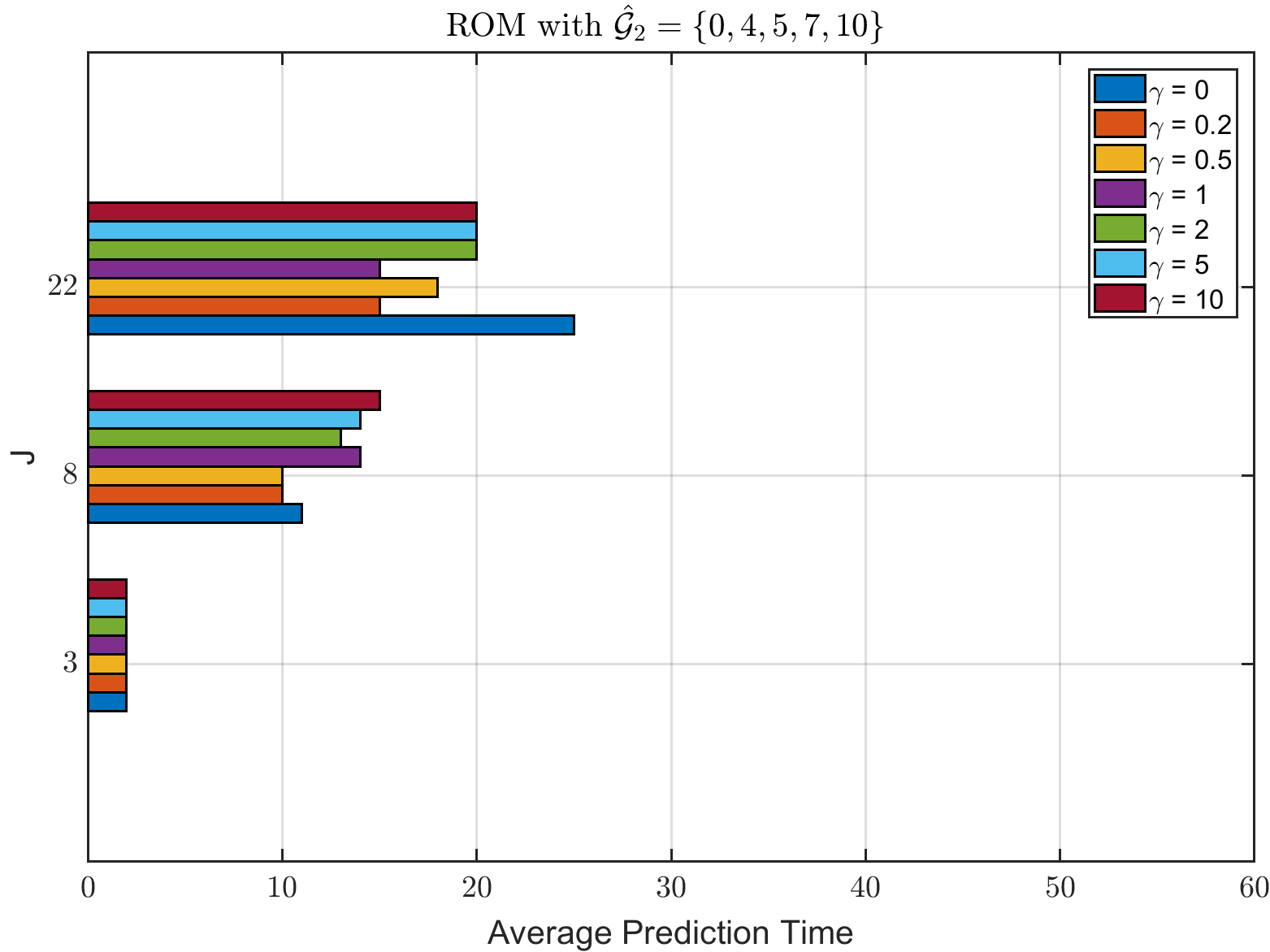}}
\vspace{0.5cm}
\centerline{\includegraphics[width=\ww, height=\hh]{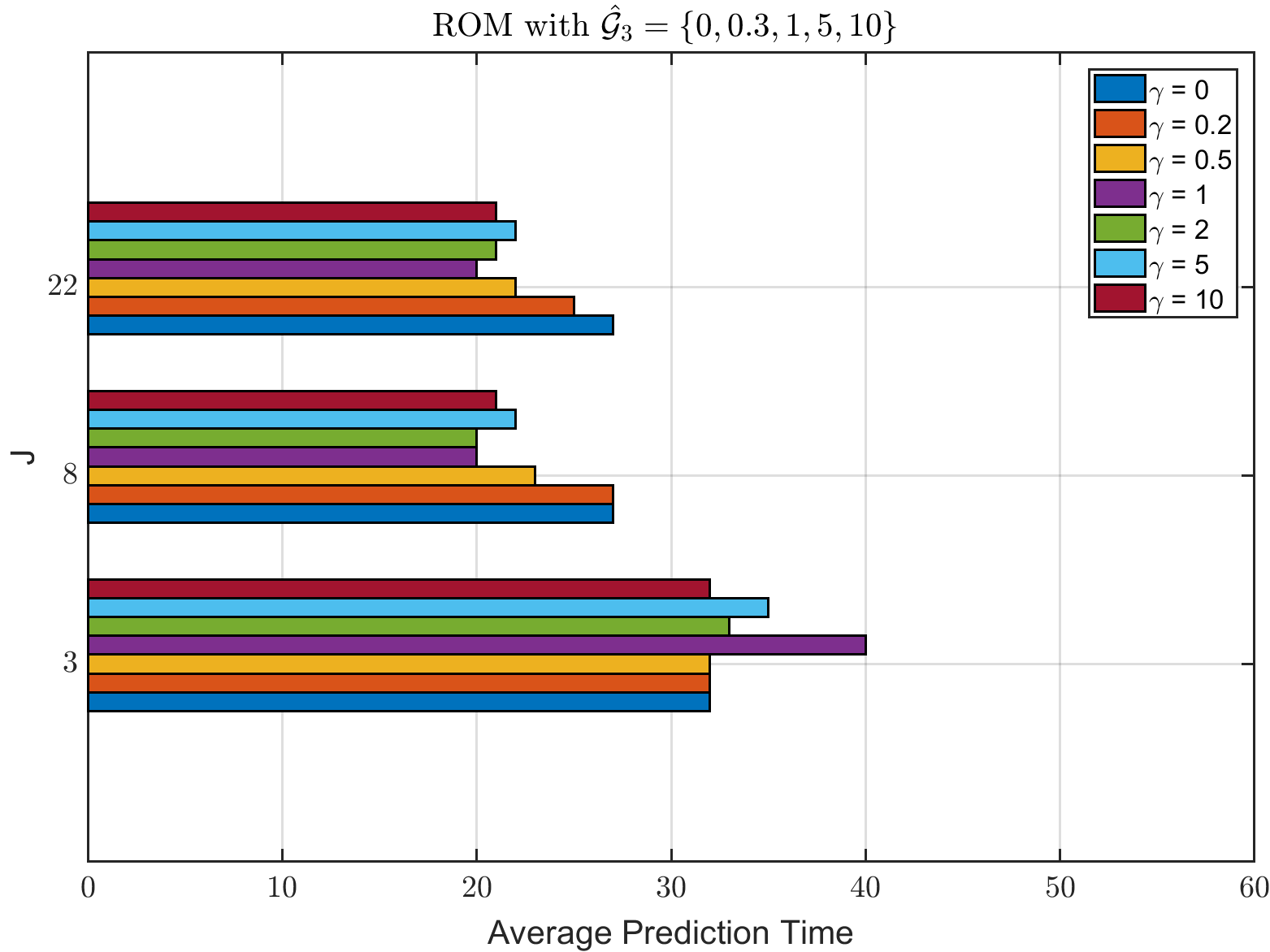}
\includegraphics[width=\ww, height=\hh]{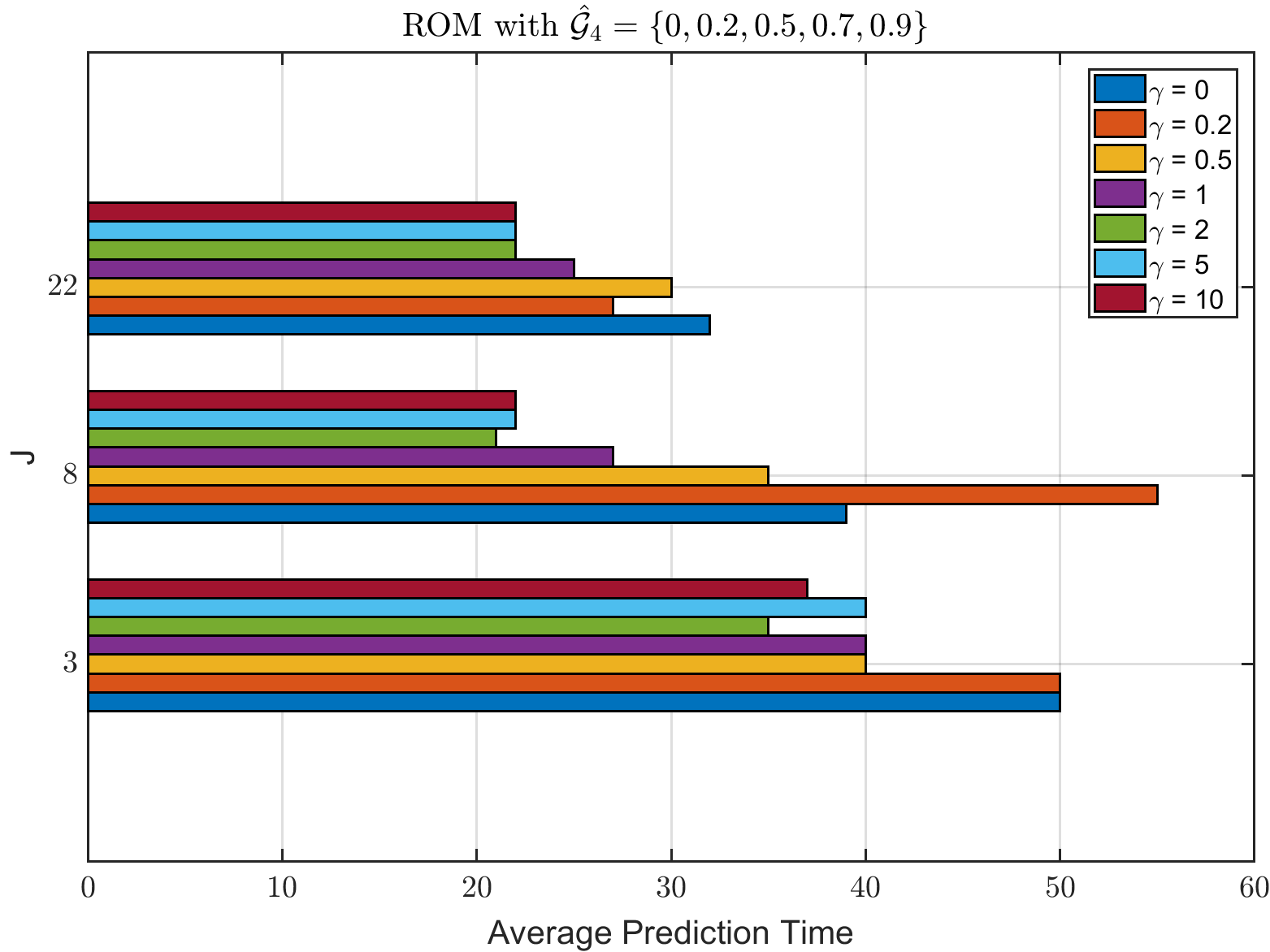}}
\caption{Averaged prediction times for trajectories of the gKS equation with different values of $\gamma$ using multi-valued ROMs constructed using gKS data with $K=5$ values of $\gamma$ and $W=1$. For each $\gamma$, we consider 3 sets of different initial conditions with $J=3,8,22$, where $J$ is the number of non-zero Fourier wavenumbers at time $t=0$ (see eq. \eqref{eq:initial_condition}).}
\label{fig3}
\end{figure} 

In particular, \rev{for all ROMs} we consider initial conditions given by \eqref{eq:initial_condition} with three different values of $J = 3, 8, 22$. Averaged prediction times for two single-parameter ROMs are depicted in Figure~\ref{fig2} and for multi-parameters ROMs in Figure~\ref{fig3}. We considered three random initial conditions for each $J$ and $\gamma$ and calculated the prediction time according to \eqref{eq:predT}. The averaged prediction time is then computed by averaging over the three initial conditions.

We observe that the ROM constructed from $\gamma=0$ data has the best performance for the chaotic regimes $\gamma=0$ and $\gamma=0.2$.  
To interpret the results for the ROM constructed from $\gamma=5$, we note that all solutions of the gKS equation with $\gamma=5$ become quasi-periodic in the long term. However, there is a short transient regime ($t \lesssim 5$) during which all solutions exhibit chaotic behavior. Therefore, as $W$ increases, the ratio of these ``chaotic'' snapshots increases in $\Phi_{\rm pod}$ constructed from an FOM applied to the gKS equation with $\gamma=5$.  

Figure~\ref{fig2} also demonstrates that as the number of ``chaotic'' snapshots in $\Phi_{\rm pod}$ increases, ROMs perform better in terms of short-term prediction across all regimes $0 \leq \gamma \leq 10$. Similarly, Figure~\ref{fig3} indicates that ROMs constructed with more ``chaotic'' snapshots perform better. In particular, in Figure~\ref{fig3}, we consider multi-value ROMs, with training parameter sets ${\hcA}_i$ chosen such that the number of ``chaotic'' snapshots increases with the index $i$.

\subsubsection{Chaotic regime and statistical behavior.}\label{s:Ch}
The gKS equation converges to the integrable KdV equation as $\gamma\to\infty$ and, thus, 
the long-term behavior of solutions becomes quasi-periodic as $\gamma$ increases (e.g. \cite{chang1993laminarizing,pradas2011rigorous,tseluiko2014weak}).
It is possible to identify three different regimes for the long-term behavior of solutions of the gKS equation. Parameter values $\gamma\in[0,0.19]$ correspond to the chaotic regime, $\gamma \gtrsim 1$ corresponds to quasi-periodic behavior, and $\gamma\in[0.2,0.9]$ is an intermediate regime where solutions can 
exhibit transient chaotic behavior for a long time. This behavior of the gKS equation was analyzed extensively numerically in, e.g., \cite{gotoda2015}.

\begin{figure}
\centerline{\includegraphics[width=\ww]{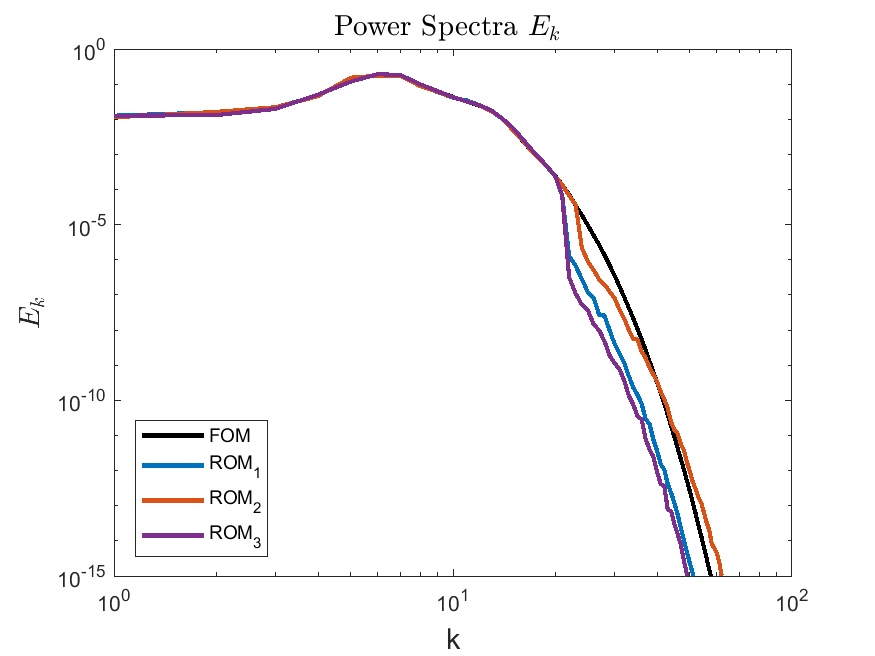}
\includegraphics[width=\ww]{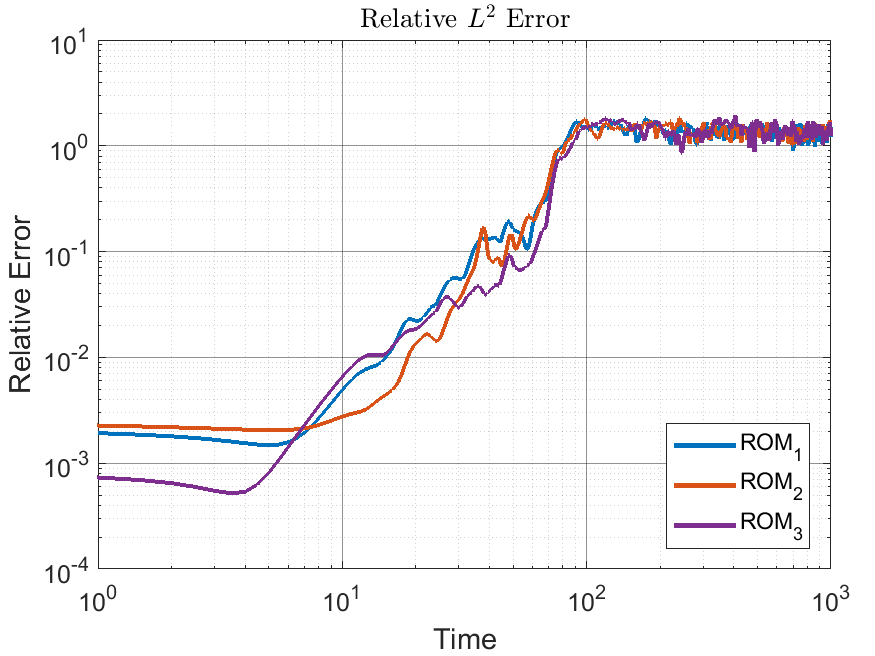}}
\caption{Power Spectra (left) and a typical $L^2$ error for a single trajectory (right) for $\gamma = 0.1$.}
\label{fig4}
\end{figure}

\begin{figure}
\centerline{\includegraphics[width=\ww]{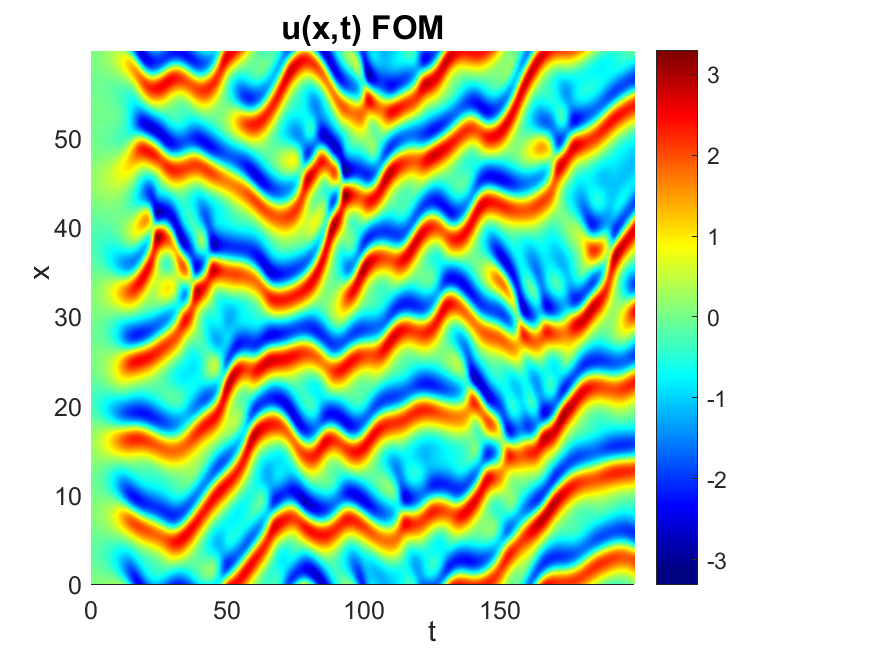}
\includegraphics[width=\ww]{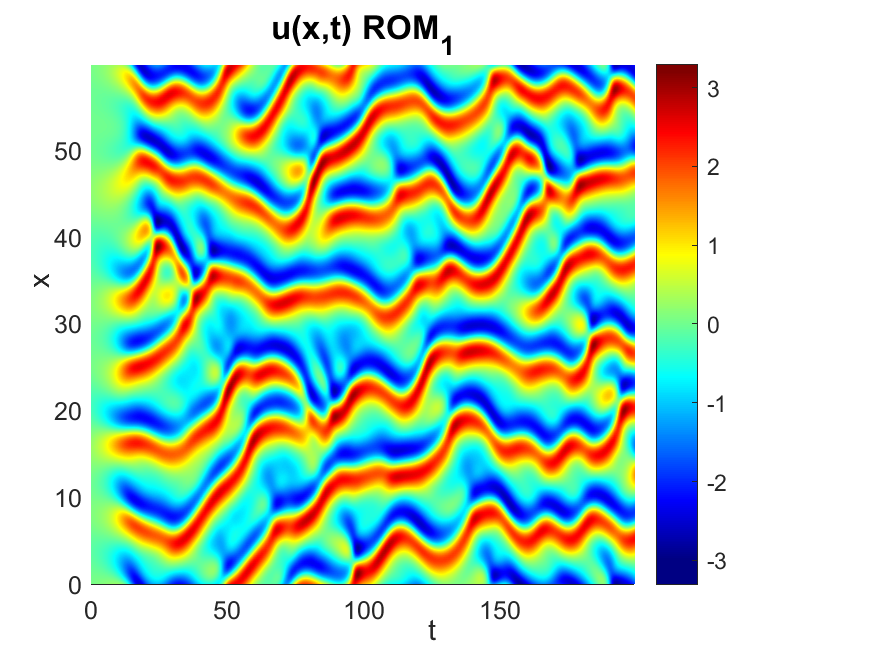}}
\centerline{\includegraphics[width=\ww]{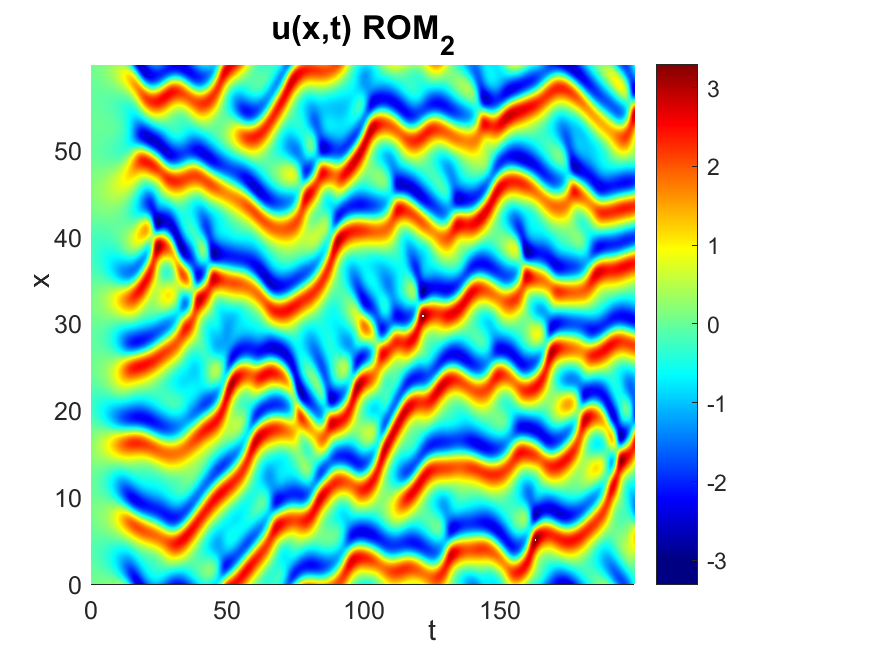}
\includegraphics[width=\ww]{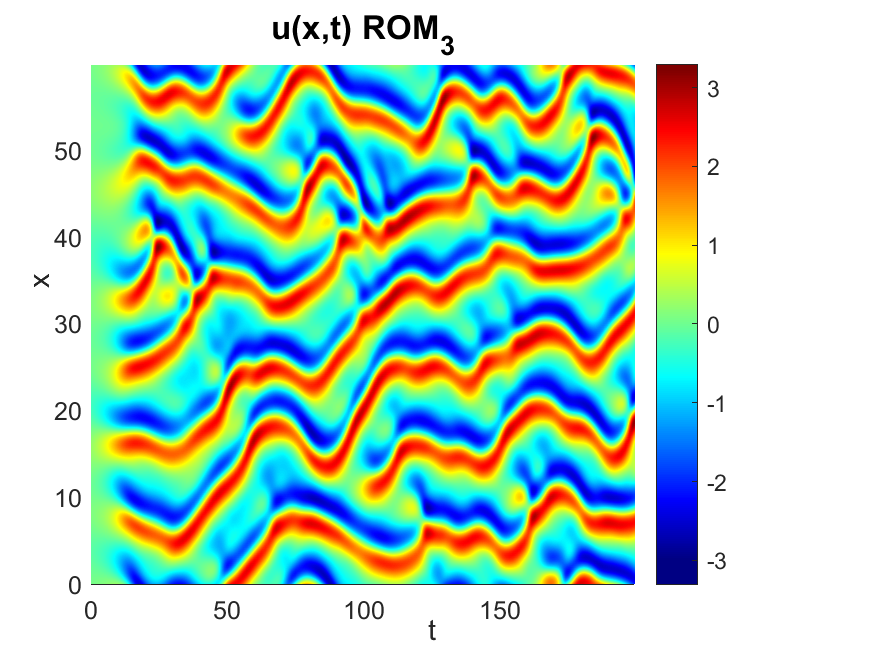}}
\caption{Comparison of an individual solution for FOM and $\text{ROM}_j$, $j=1,2,3$, for $\gamma = 0.1$.}
\label{fig5}
\end{figure}

Although, as shown previously, ROM accuracy deteriorates for individual trajectories at times around $t=100$,
we can compare the long-time statistical behavior of ROMs for the gKS equation. To this end, we compute the power spectra 
over a long time. The power spectra are given by
\[
E_k = \frac1T \int\limits_0^T |\hat{u}_k(t)|^2 \, dt,
\]
where $\hat{u}_k(t)$ are the coefficients of the Discrete Fourier Transform of the solution $u(x,t)$.

\medskip
\textbf{Definition of $\text{ROM}_j$ with $j=1,2,3$.}
 We consider three particular ROMs discussed previously:\\[0.3ex] 
 $\text{ROM}_1$ is constructed from snapshots computed with $\gamma=0$, $W=1$, and $K=1$; \\
 $\text{ROM}_2$ is constructed using $\gamma=5$, $W=250$, and $K=1$;  \\
 $\text{ROM}_3$ is multi-valued constructed using $K=5$ and ${\hcA}_4 = \{0, 0.2, 0.5, 0.7, 0.9\}$.
\smallskip

Figure~\ref{fig4} depicts the power spectra for the FOM and $\text{ROM}_j$ ($j=1,2,3$), as well as the $L^2$ error of the individual solution between the FOM and $\text{ROM}_j$ for $\gamma=0.1$. Relative $L^2$ errors are shown on a log-log scale to emphasize the behavior near $t=0$ and suppress noisy oscillations in the graph of the relative error.  
Figure~\ref{fig5} presents a comparison of one particular trajectory on the time interval $[0,200]$. Numerical solutions computed using the FOM and $\text{ROM}_j$ appear very similar for $t \in [0,80]$. This observation is supported by the Averaged Relative Error depicted in Figure~\ref{fig4}.  

Larger differences between the numerical solutions computed using the FOM and $\text{ROM}_j$ become clearly visible on the time interval $t \in [100,200]$. However, the statistical behavior remains very similar. Consequently, the power spectra for the FOM and $\text{ROM}_j$ (Figure~\ref{fig4}) are almost identical for low wavenumbers, with only a small discrepancy observed at higher wavenumbers in the $\text{ROM}_j$ simulations.  

All three ROMs successfully reproduce the solution on the time interval $[0,80]$ (short-term prediction) and correctly capture the statistical behavior of the solutions for longer times. This outcome is expected since the gKS equation exhibits chaotic solutions for $\gamma=0.1$, and numerical errors eventually lead to rapidly divergent trajectories of the ROM and FOM solutions. Nevertheless, all three ROMs accurately capture the behavior of the attractor of the gKS equation.

\subsubsection{Transient regime and persistent patterns.}\label{s:Tr}
For $\gamma>0.2$, statistical properties of the gKS equation become 
dependent on the initial condition. Therefore, we do not present the averaged power spectra for values of $\gamma>0.2$. For $\gamma\in[0.2,0.9]$, solutions can exhibit long chaotic-like transition periods, but eventually, solutions become quasi-periodic and persistent patterns emerge.

\begin{figure}
\centerline{\includegraphics[width=\ww]{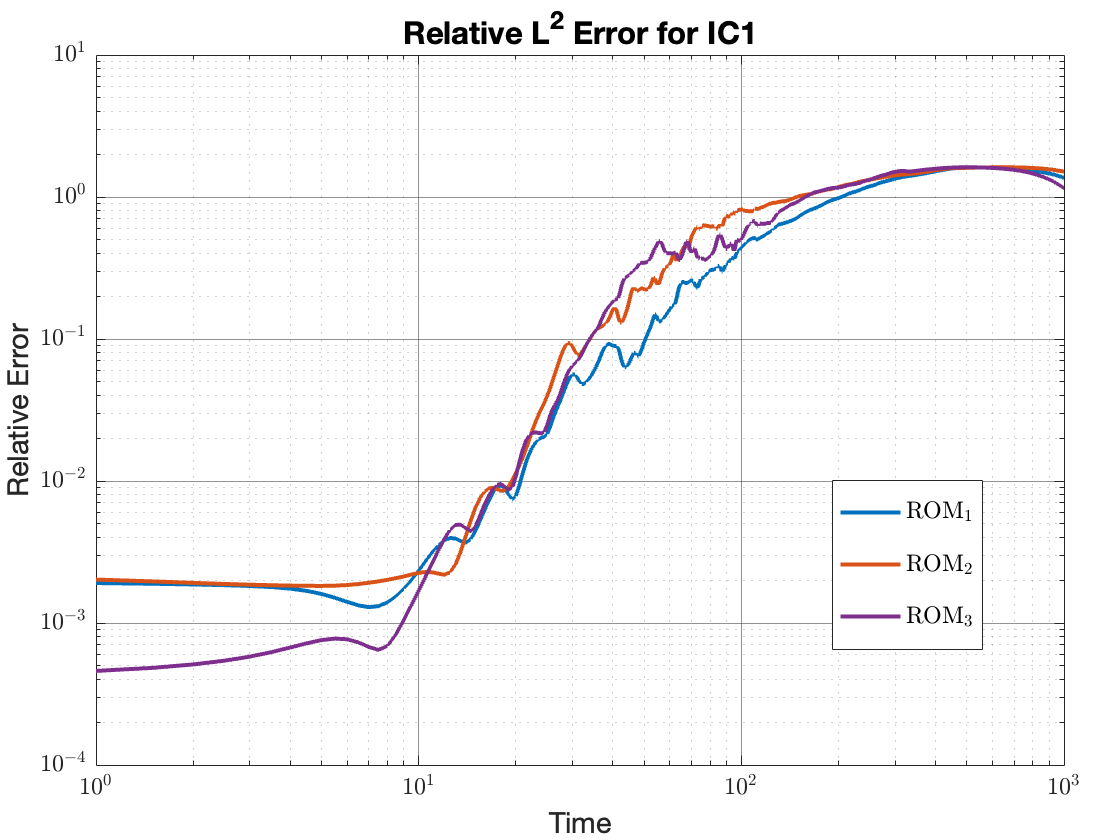}
\includegraphics[width=\ww]{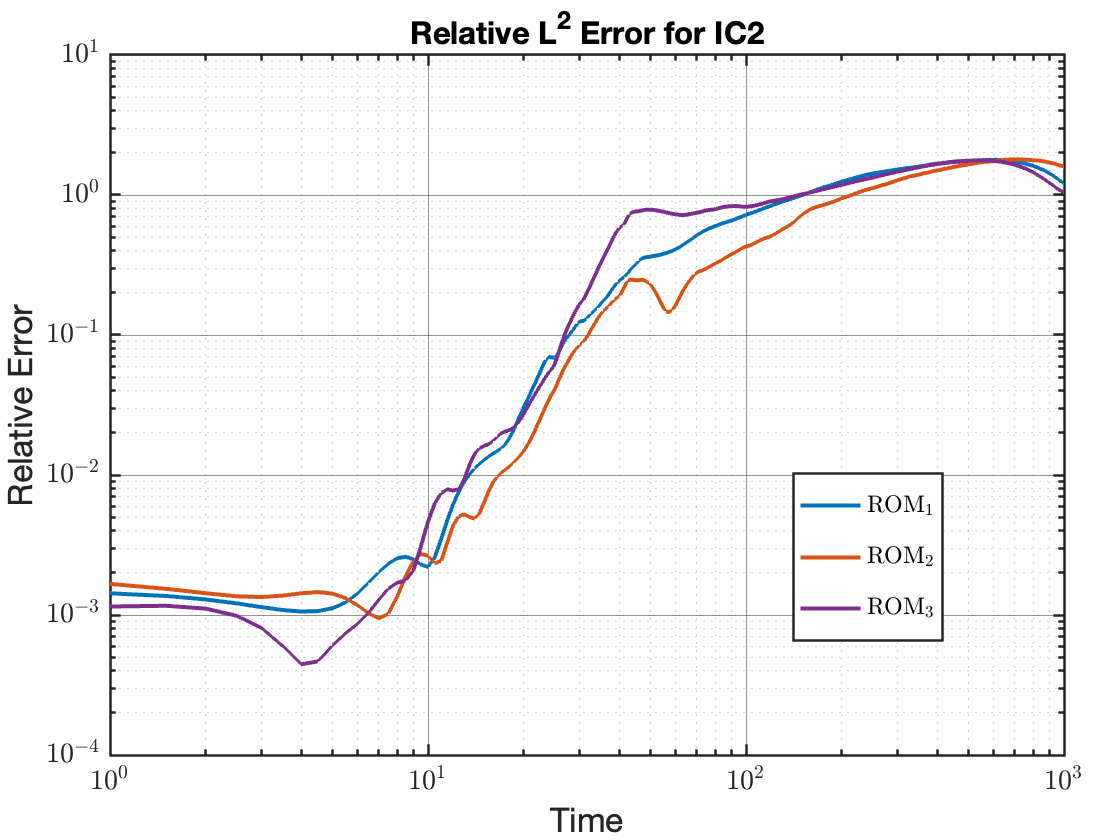}}
\caption{$L^2$ error for two particular trajectories (given by different initial conditions; left: IC1, right: IC2) in simulations of FOM and $\text{ROM}_j$, $j=1,2,3$ with $\gamma = 0.7$.}
\label{fig6}
\end{figure}

\begin{figure}
\centerline{\includegraphics[width=\ww]{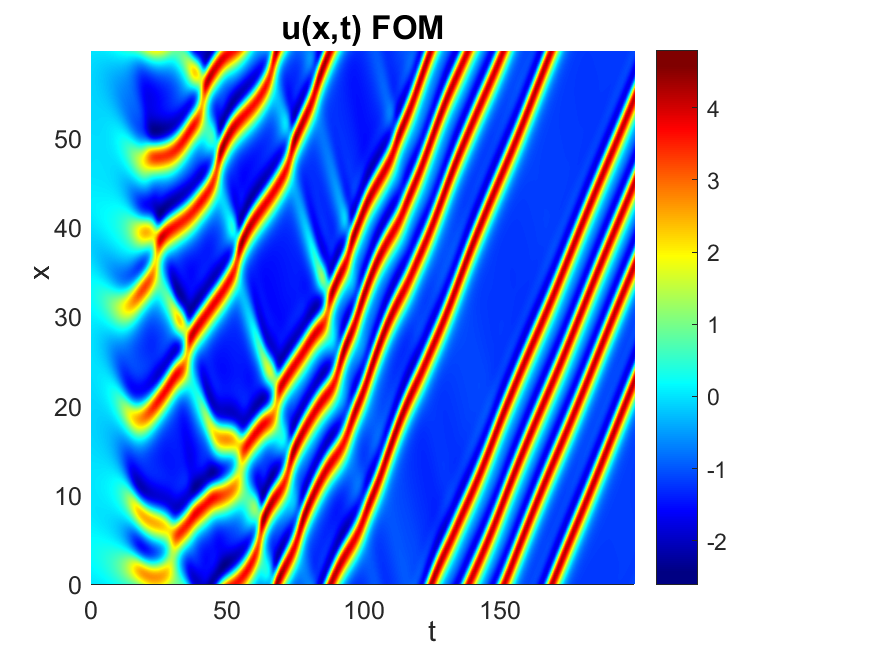}
\includegraphics[width=\ww]{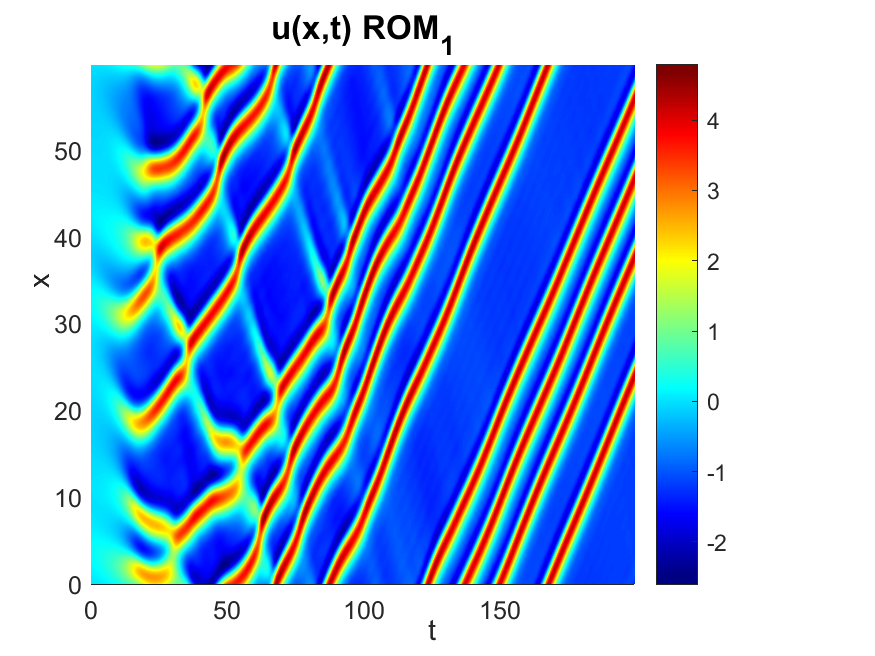}}
\centerline{\includegraphics[width=\ww]{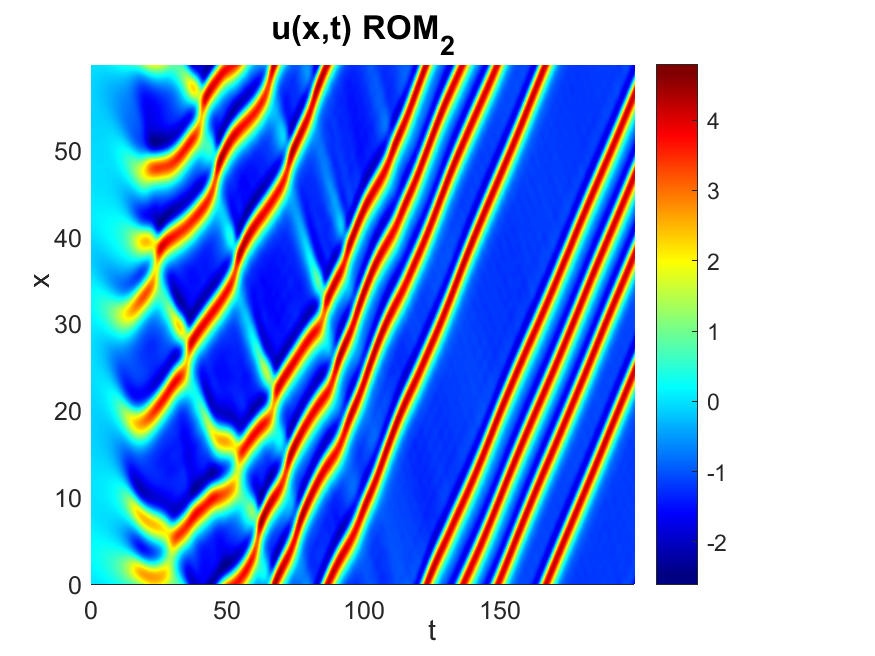}
\includegraphics[width=\ww]{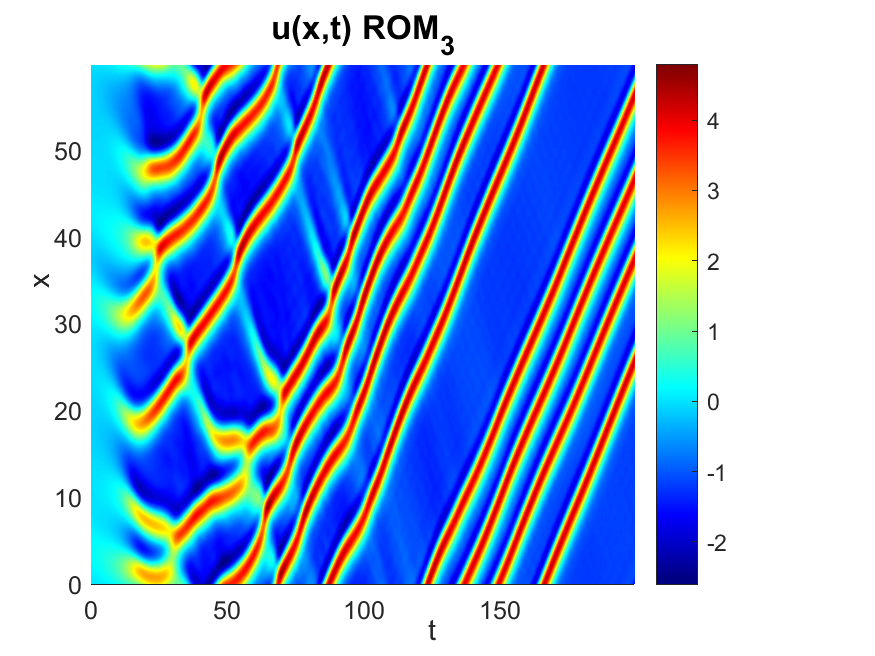}}
\caption{Snapshots fo solutions in simulations of FOM and $\text{ROM}_j$, $j=1,2,3$ with $\gamma=0.7$. Initial condition IC1.}
\label{fig7}
\end{figure}
\begin{figure}
\centerline{\includegraphics[width=\ww]{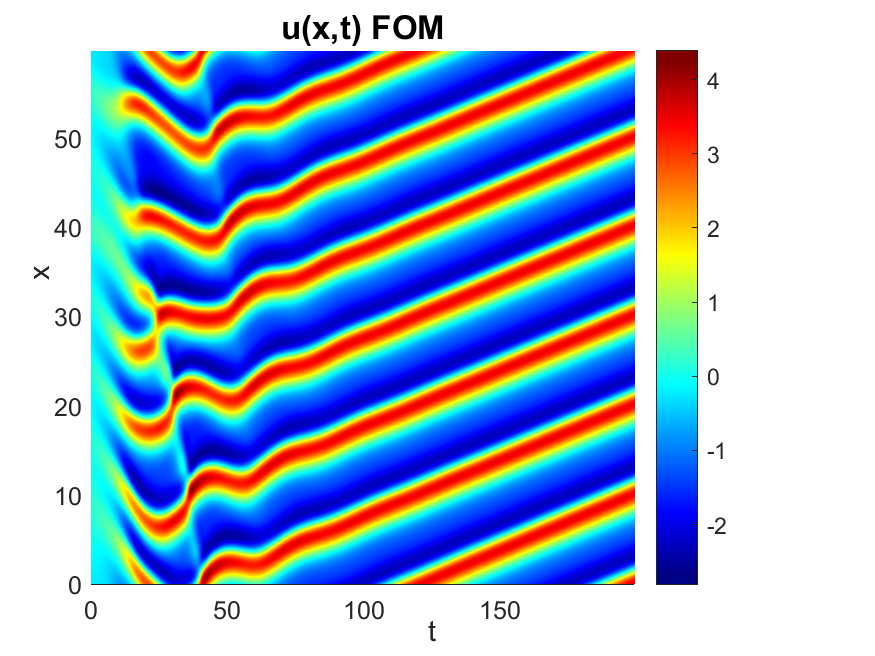}
\includegraphics[width=\ww]{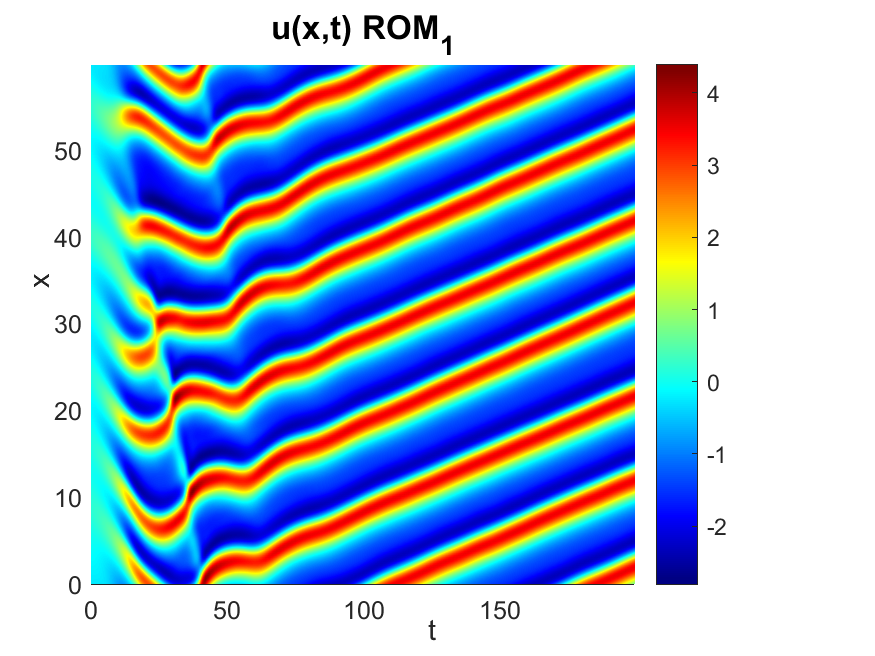}}
\centerline{\includegraphics[width=\ww]{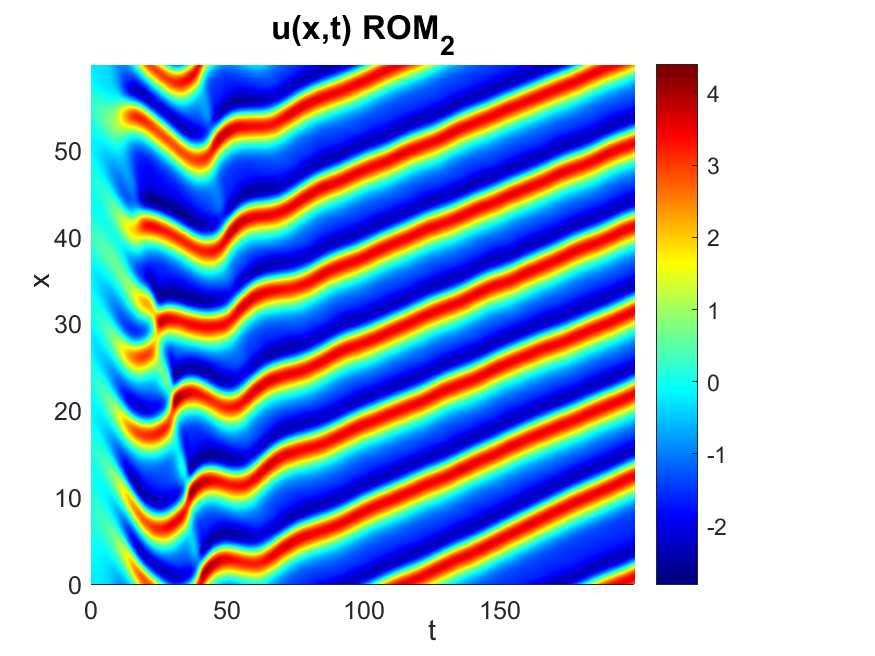}
\includegraphics[width=\ww]{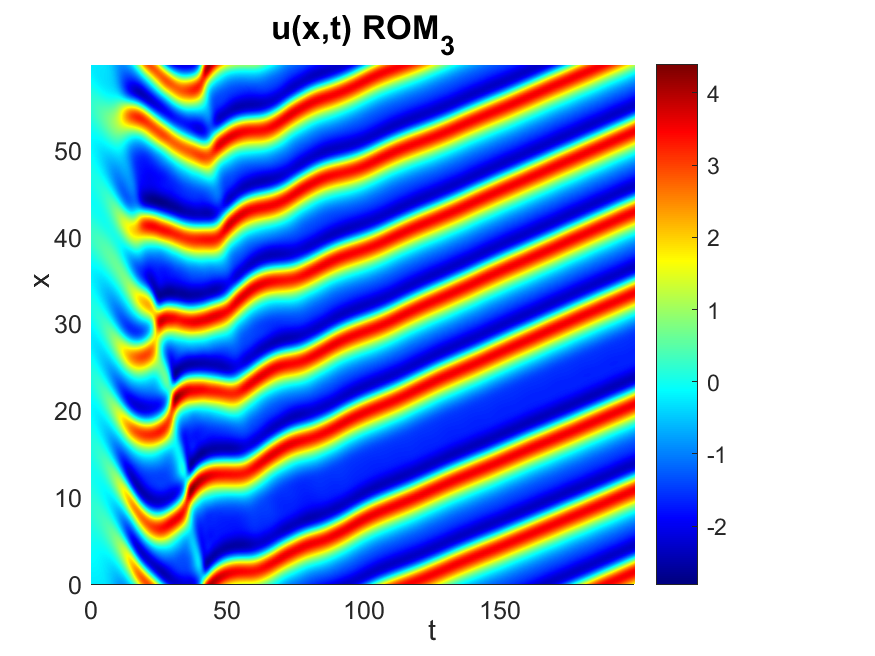}}
\caption{Snapshots fo solutions in simulations of FOM and $\text{ROM}_j$, $j=1,2,3$ with $\gamma=0.7$. Initial condition IC2.}
\label{fig8}
\end{figure}

We present the results of ROM and FOM numerical simulations for $\gamma=0.7$ and two particular initial conditions (IC1 and IC2) in Figures~\ref{fig6}, \ref{fig7}, and \ref{fig8}. For both initial conditions, the solutions demonstrate persistent traveling wave-like structures that form beyond the transient stage at longer times. The speed, width, and amplitude of these traveling structures are clearly dependent on the initial conditions.

We observe that all three ROMs accurately reproduce the solution up to $t \simeq 10$, after which the relative $L^2$ error between the ROM and FOM solutions begins to grow, reaching $O(1)$ around $t=100$; see Figure~\ref{fig6}. Nevertheless, the persistent solution patterns are very well reproduced by the ROMs when the long-term behavior of solutions is of interest. The larger errors for $t > 100$ can be attributed to small phase shifts in the traveling-wave solutions between the ROMs and the FOMs.

We conclude that the error norm, a commonly used statistic for assessing the quality of a ROM, may be inadequate for systems that exhibit laminar solutions with persistent patterns emerging after a transition stage.

\subsubsection{Laminar regime and quasi-periodic solutions.} \label{s:La}

For larger $\gamma$ the gKS system becomes deterministic with solutions exhibiting a short transition period before reaching a quasi-periodic state. 
Both the transition period and the quasi-periodic state are initial condition dependent. The later makes it hard to predict long term solution behavior especially  if the initial condition is not within (or very close to) the training set of initial conditions. This phenomenon is illustrated by the following example of recovering gKS solution for $\gamma=3$. 

\begin{figure}
\centerline{\includegraphics[width=\ww]{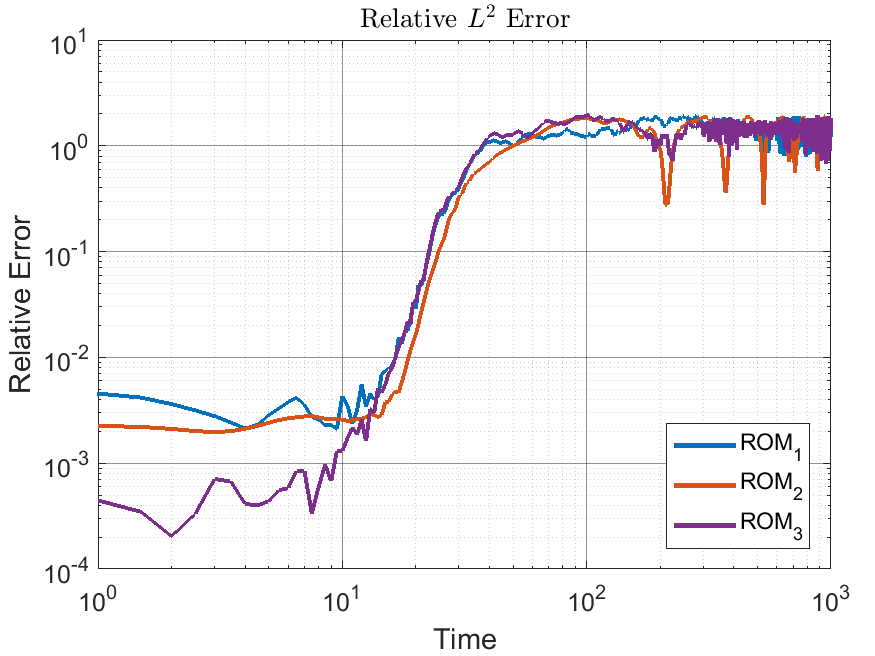}
\includegraphics[width=\ww]{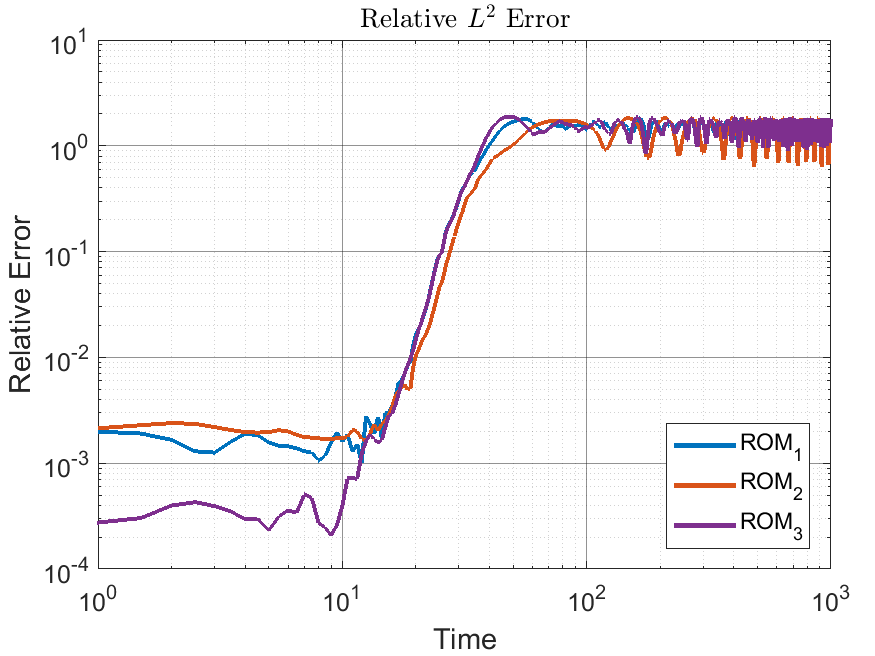}}
\caption{$L^2$ error for two particular trajectories (given by different initial conditions; left: IC1, right: IC2) in simulations of FOM and $\text{ROM}_j$, $j=1,2,3$ with $\gamma = 3$.}
\label{fig9}
\end{figure}
\begin{figure}
\centerline{\includegraphics[width=\ww]{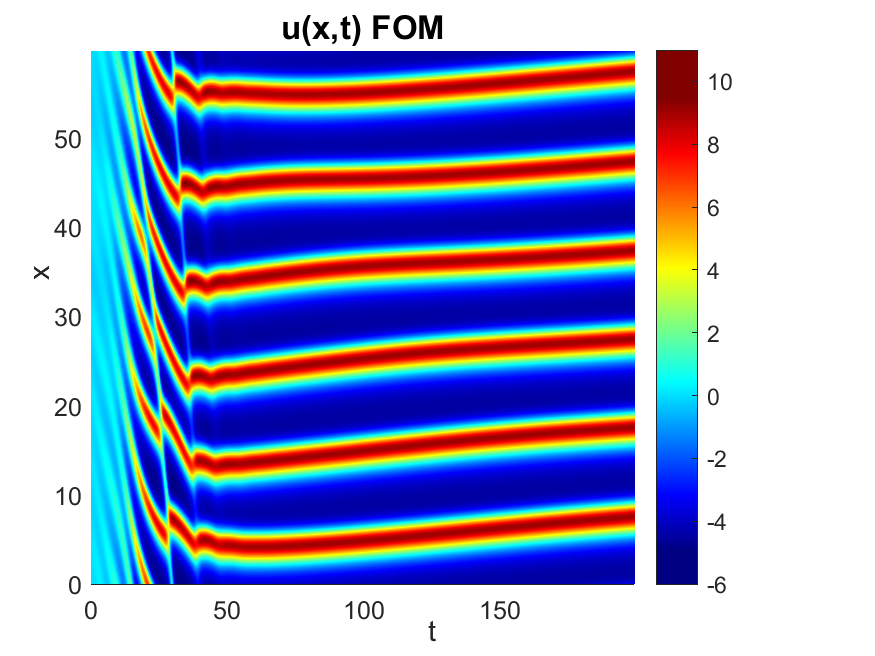}
\includegraphics[width=\ww]{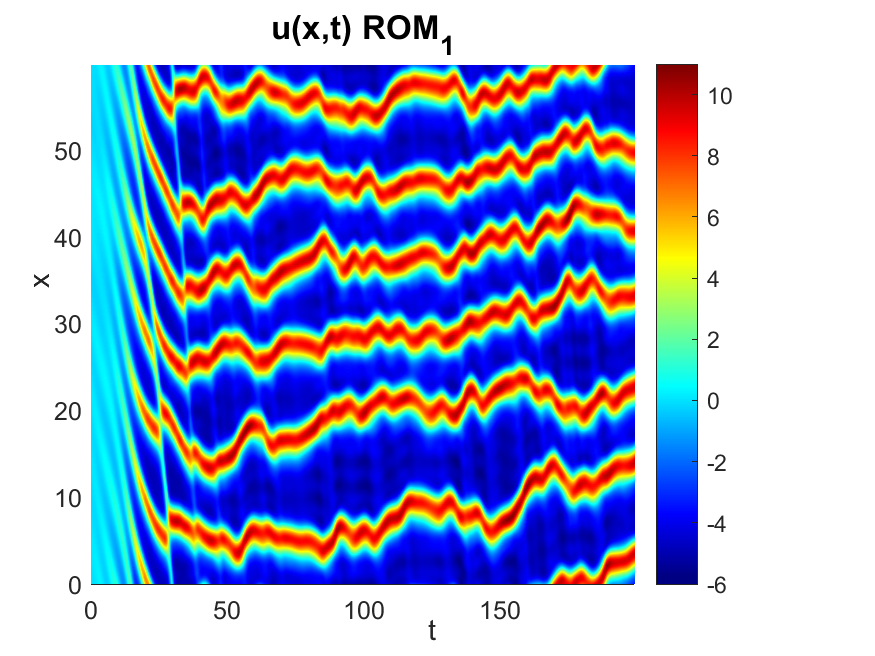}}
\centerline{\includegraphics[width=\ww]{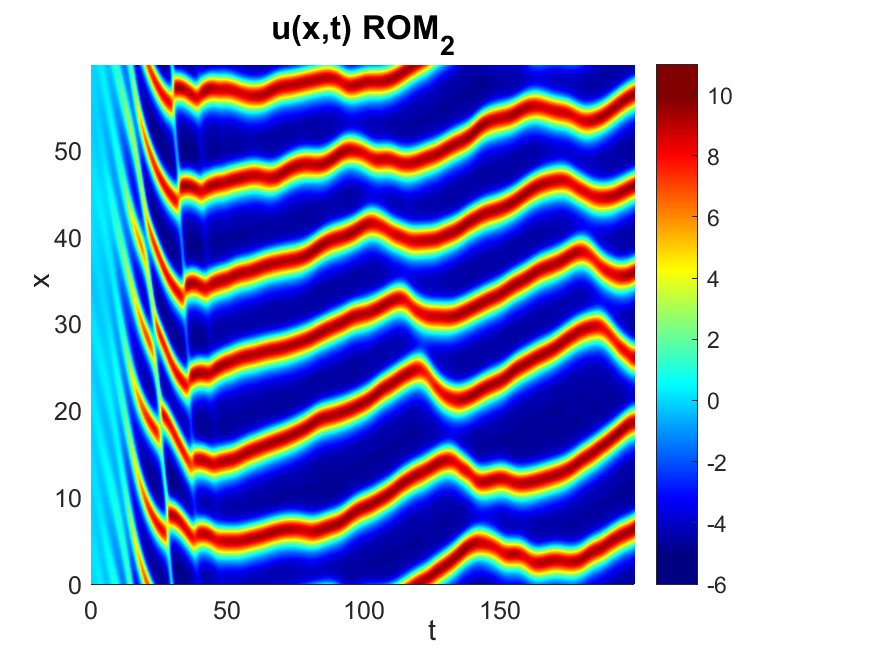}
\includegraphics[width=\ww]{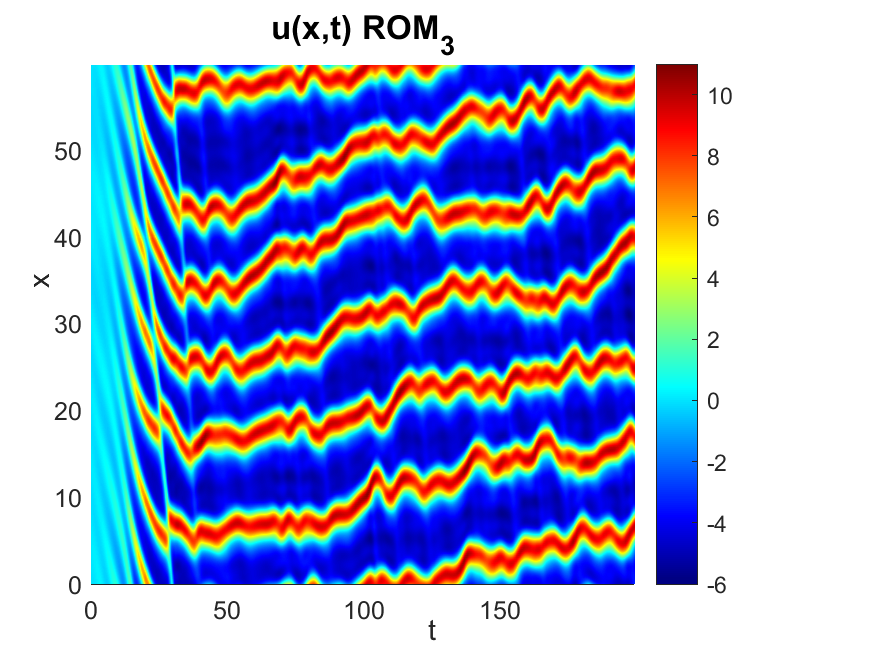}}
\caption{Snapshots fo solutions in simulations of FOM and $\text{ROM}_j$, $j=1,2,3$ with $\gamma=3$. Initial condition IC1.
$\text{ROM}_3$ is constructed with ${\hcA}_4 = \{0, 0.2, 0.5, 0.7, 0.9\}$.}
\label{fig10}
\end{figure}
\begin{figure}
\centerline{\includegraphics[width=\ww]{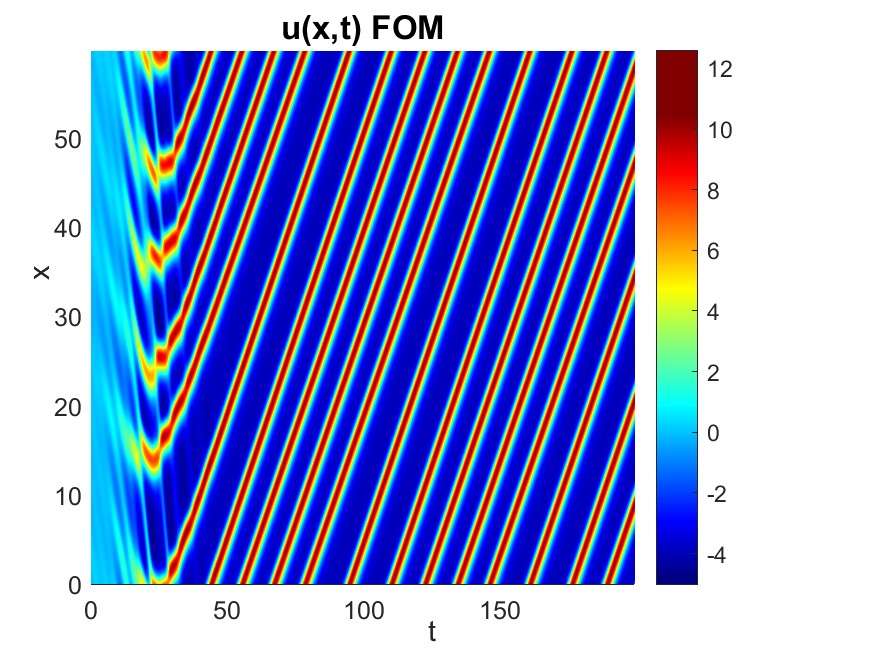}
\includegraphics[width=\ww]{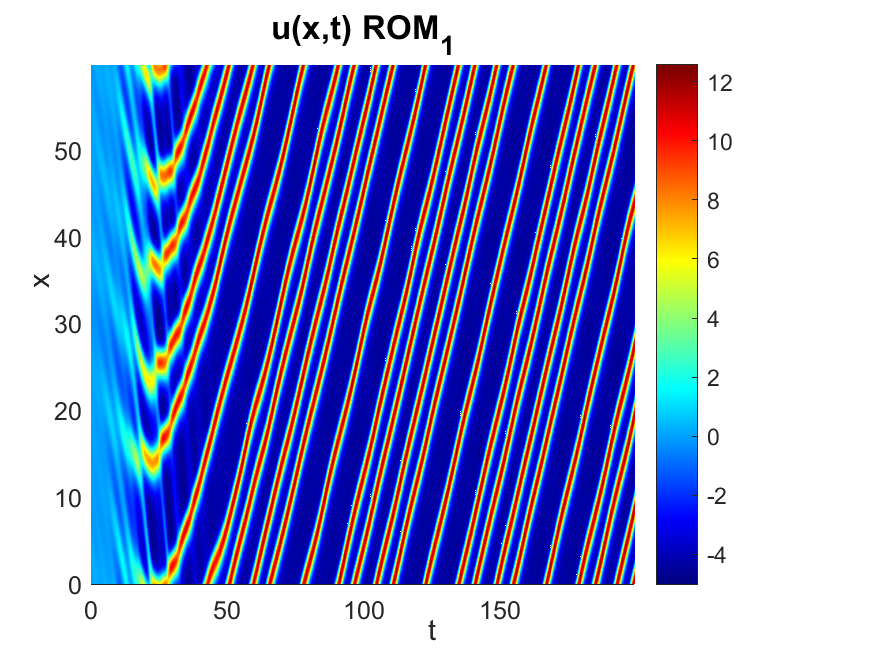}}
\centerline{\includegraphics[width=\ww]{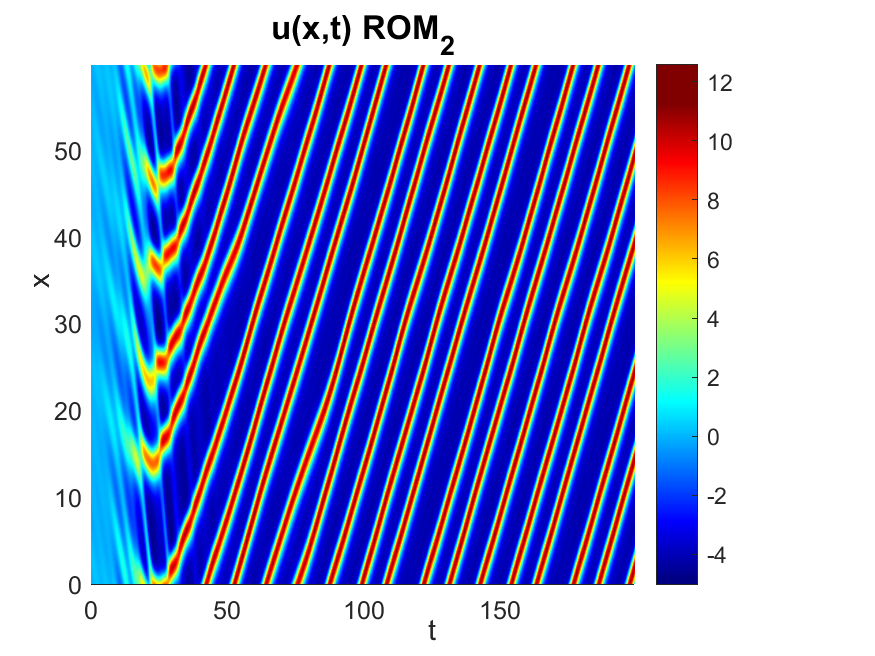}
\includegraphics[width=\ww]{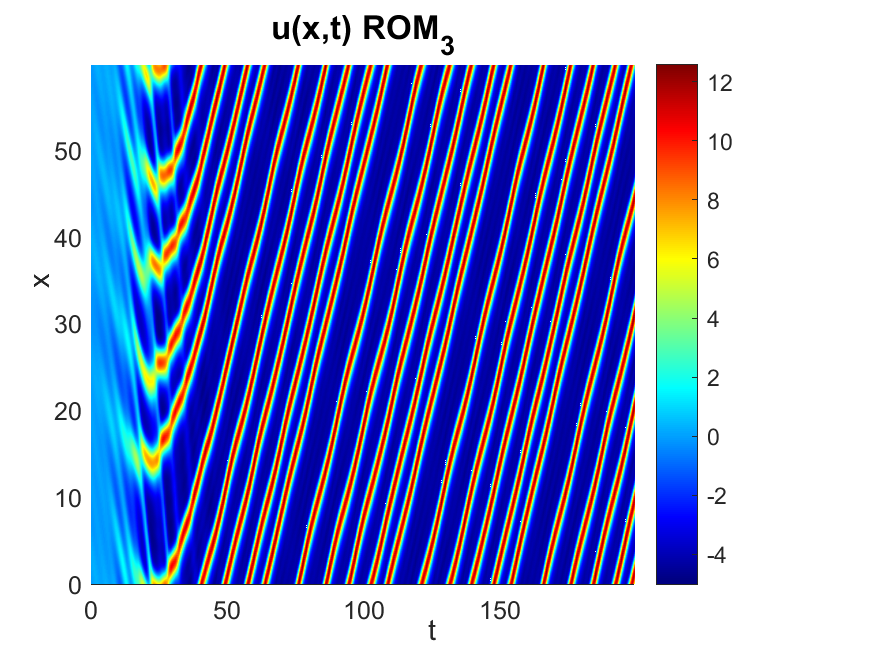}}
\caption{Snapshots fo solutions in simulations of FOM and $\text{ROM}_j$, $j=1,2,3$ with $\gamma=3$. Initial condition IC2. 
$\text{ROM}_3$ is constructed with ${\hcA}_4 = \{0, 0.2, 0.5, 0.7, 0.9\}$.}
\label{fig11}
\end{figure}

Figures~\ref{fig9}--\ref{fig10} show the results of $\text{ROM}_j$, $j=1,2,3$, and FOM numerical simulations for $\gamma=3$ and two particular initial conditions (IC1 and IC2; \rev{both not in the training set of the initial conditions}). All three ROMs accurately reproduce the solution up to $t \simeq 10$, after which the relative $L^2$ error between the ROM and FOM solutions begins to grow, reaching $O(1)$ around $t=40$; see Figure~\ref{fig9}. Unlike the transitional case with persistent solution patterns, the ROMs fail to accurately reproduce the quasi-periodic state of the solution. However, we observe that $\text{ROM}_2$ performs somewhat better in the long term, offering the best recovery of the quasi-periodic states among all ROMs considered.

We note that other variants of ROMs described at the beginning of the section were also applied in this example. In particular, variants of multi-parameter ROMs with training sets ${\hcA}_i$, $i=1,2,3$, were tested. They all demonstrated results no better than $\text{ROM}_1$ or $\text{ROM}_3$. Notably, the results for the multi-parameter ROM with training set ${\hcA}_1$ were the least accurate, despite $\gamma=3$ being among the training values for this ROM. We omit presenting these results for brevity. A similar pattern was observed when ROM and FOM solutions were compared for gKS with $\gamma=5$.

\subsection{POD-DEIM and POD ROM comparison} \label{s:DEIM}
A hyper-reduction technique is generally required for effective handling of non-linear terms in a projection-based ROM. In this final set of experiments, we demonstrate that applying the discrete empirical interpolation method (DEIM) leads to a reduced model that is as accurate as the (computationally expensive) true projection. In other words, we compare ROMs given by \eqref{eqn:genericROM} and \eqref{POD-D-ROM}.

% Case Study: gamma = 0
\begin{figure}[H]
    \centering
    \begin{subfigure}{0.4\textwidth}
        \centering
        \includegraphics[width=\textwidth]{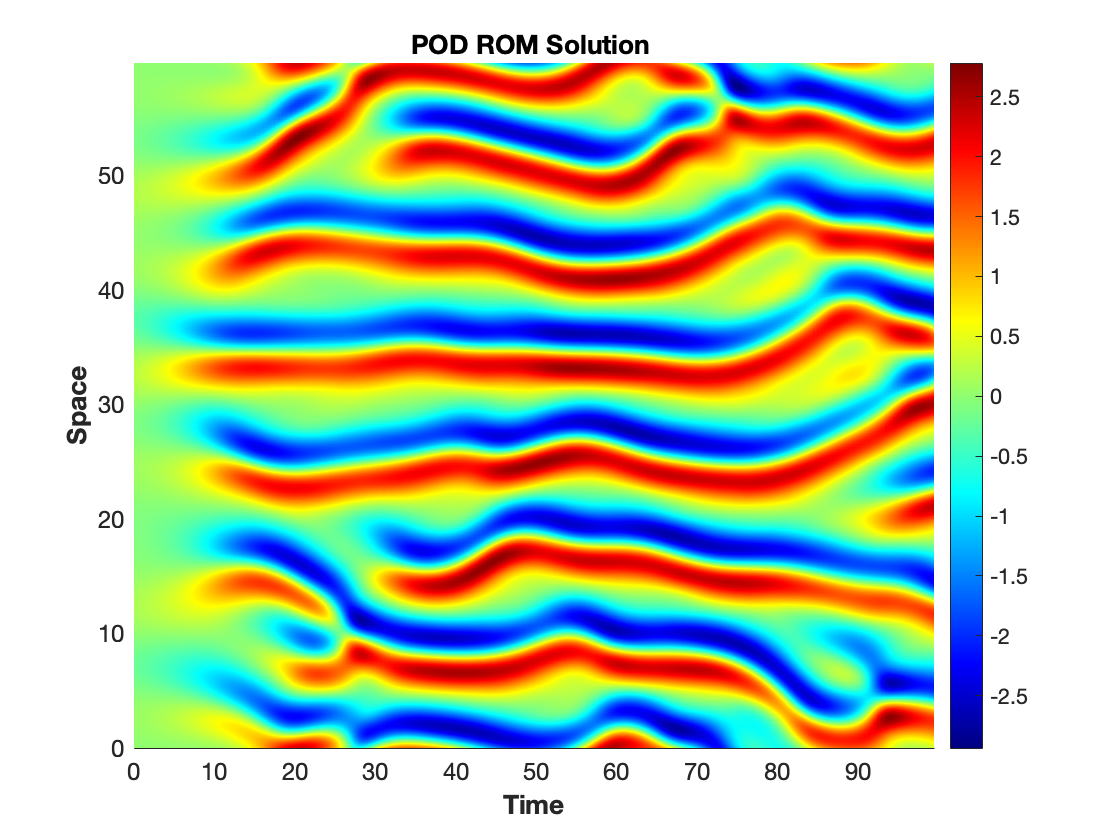}
        \caption{POD ROM Solution }
        \label{fig:fig1}
    \end{subfigure}
    \hskip3ex
    \begin{subfigure}{0.4\textwidth}
        \centering
        \includegraphics[width=\textwidth]{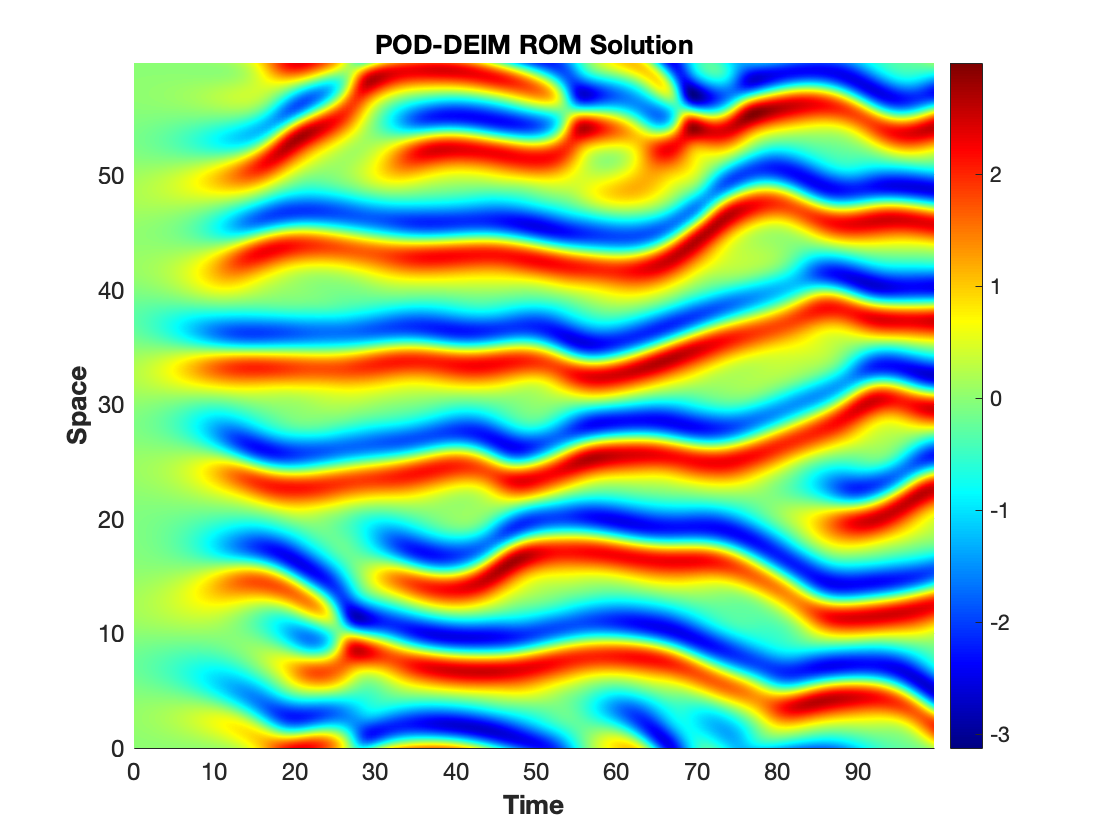}
        \caption{POD-DEIM ROM Solution}
        \label{fig:fig2}
    \end{subfigure}
    \caption{POD-DEIM ROM vs POD ROM for $\gamma=0$.}
    \label{fig:all_figures_gamma_0}
\end{figure}

% Case Study: gamma = 0.5
\begin{figure}[H]
    \centering
    \begin{subfigure}{0.4\textwidth}
        \centering
        \includegraphics[width=\textwidth]{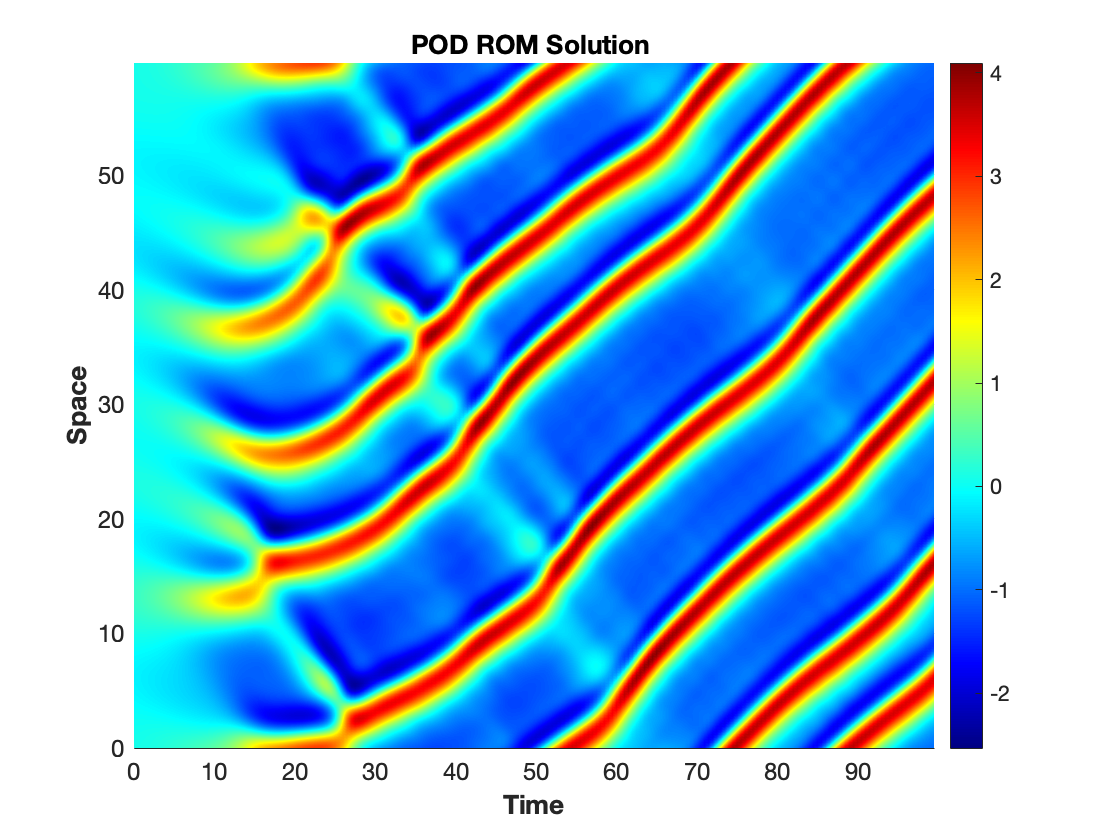}
        \caption{POD ROM Solution}
        \label{fig:fig1_0.5}
    \end{subfigure}
  \hskip3ex
    \begin{subfigure}{0.4\textwidth}
        \centering
        \includegraphics[width=\textwidth]{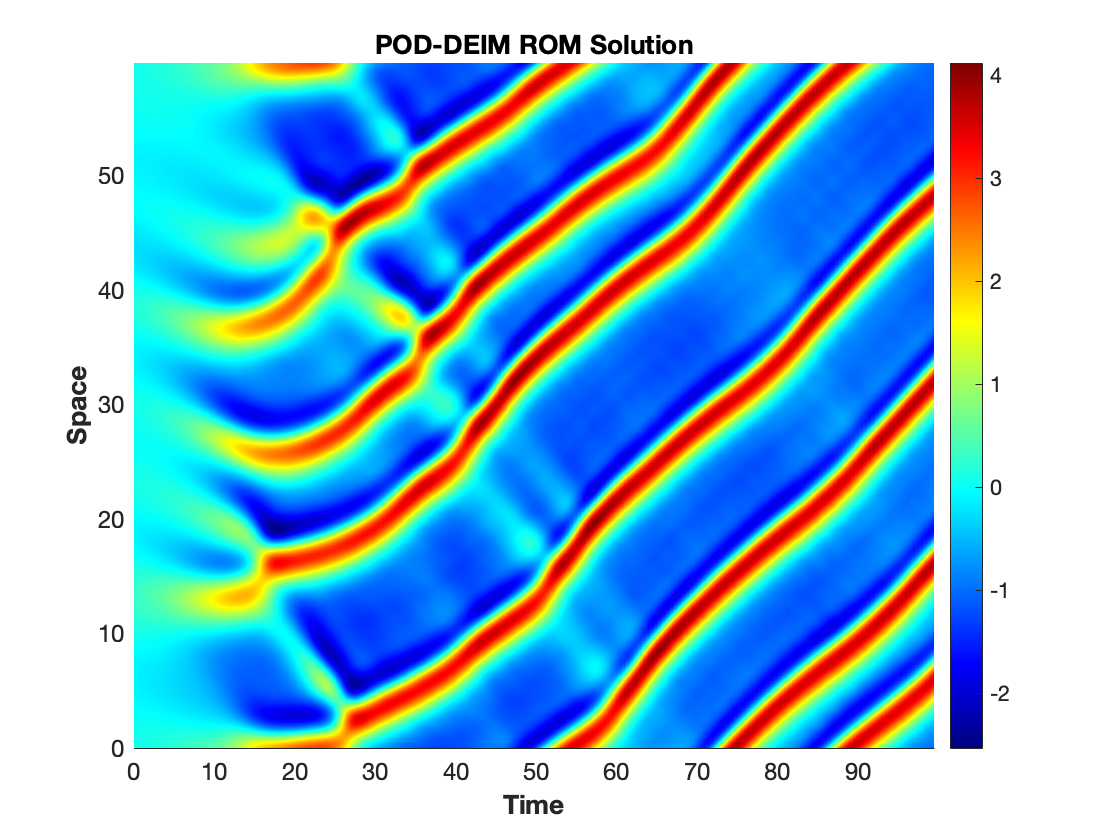}
        \caption{POD DEIM ROM Solution}
        \label{fig:fig2_0.5}
    \end{subfigure}
    \caption{POD-DEIM ROM vs POD ROM  for $\gamma=0.5$ .}
    \label{fig:all_figures_gamma_0.5}
\end{figure}

% Case Study: gamma = 2
\begin{figure}[H]
    \centering
    \begin{subfigure}{0.4\textwidth}
        \centering
        \includegraphics[width=\textwidth]{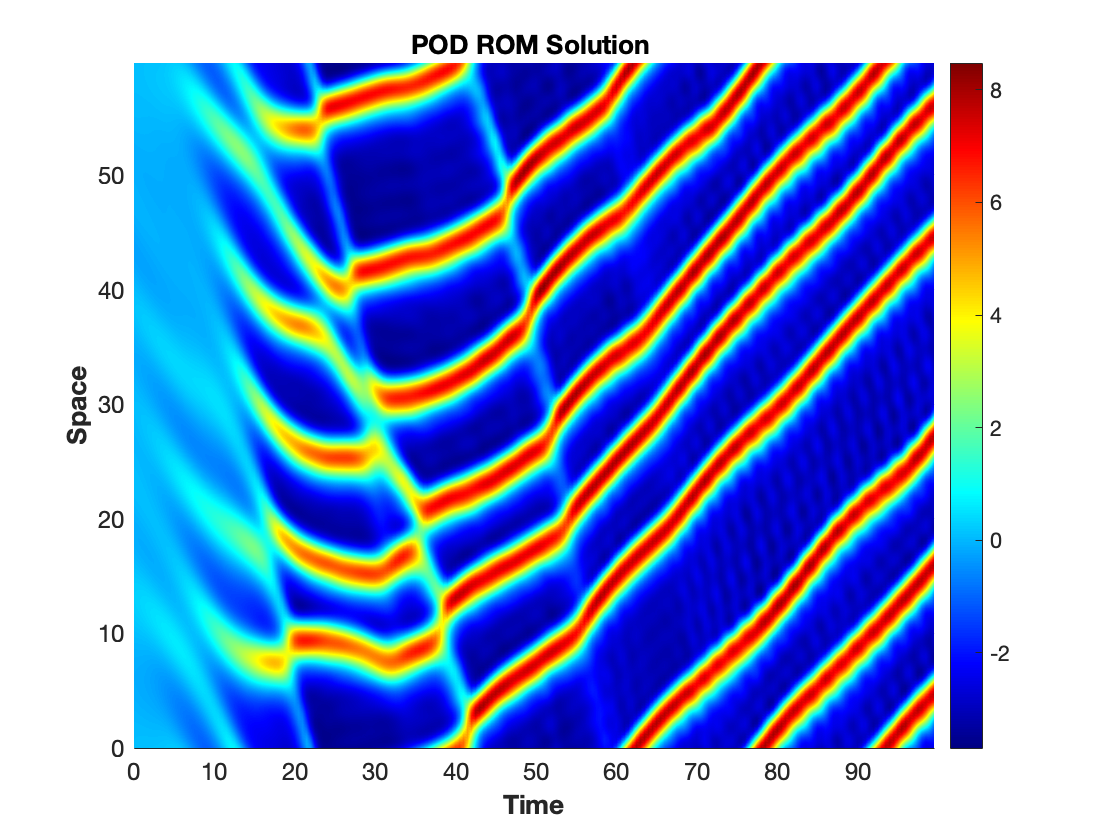}
        \caption{POD ROM Solution}
        \label{fig:fig1_02}
    \end{subfigure}
    \hskip3ex
    \begin{subfigure}{0.4\textwidth}
        \centering
        \includegraphics[width=\textwidth]{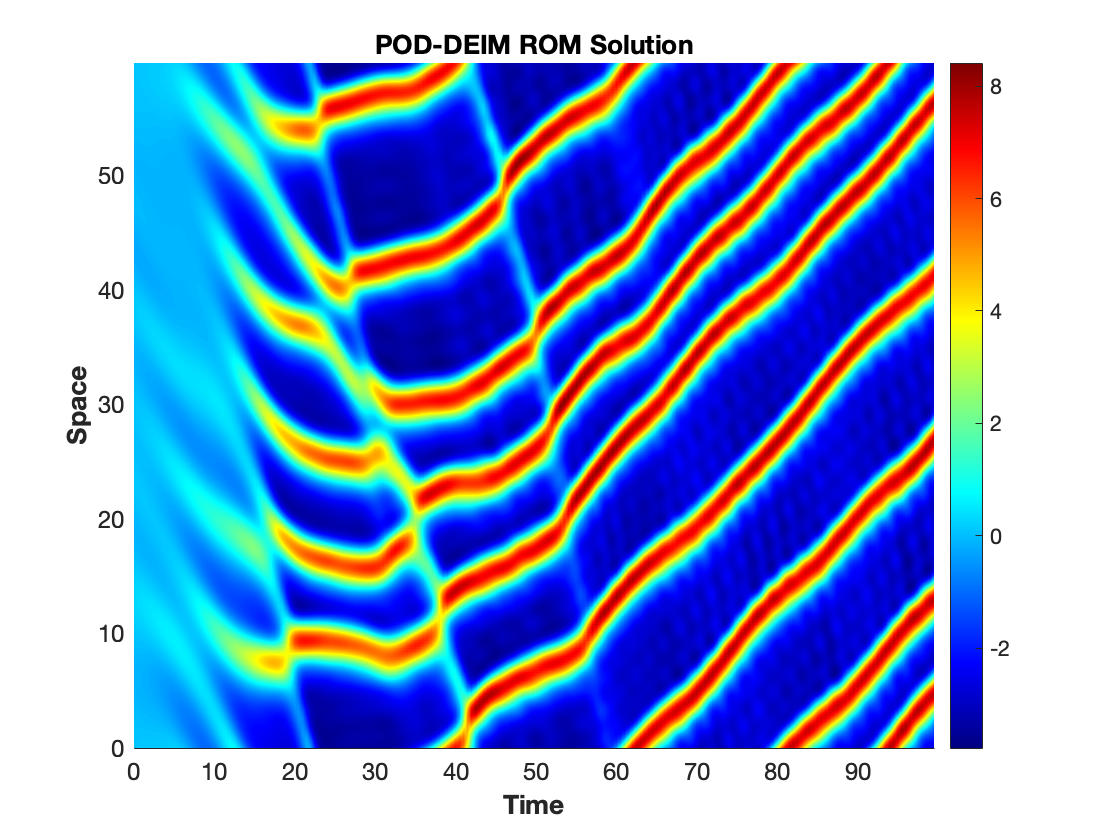}
        \caption{POD DEIM ROM Solution}
        \label{fig:fig2_02}
    \end{subfigure}
    \caption{POD-DEIM ROM vs POD ROM for $\gamma= 2$ .}
    \label{fig:all_figures_gamma_02}
\end{figure}

Figures~\ref{fig:all_figures_gamma_0}--\ref{fig:all_figures_gamma_02} present results with POD--DEIM and POD ROMs applied to predict solutions of the KS and gKS equations for \(\gamma = 0\), \(0.5\), and \(2\). We show the results only for the single-valued ROM built from \(\gamma = 0\) snapshots, as the results for ROMs built from other snapshot collections were also almost indistinguishable between the POD--DEIM and POD variants.

%%%%%%%%%%%%%%%%%%%%%%%%%%%%%%%%%
\section{Conclusions}
\label{s:conc} 
For the example of the KS and gKS equations, we investigated the ability of parametric POD-ROMs to predict solutions over a range of parameters that include various regimes exhibited by the KS and gKS solutions: chaotic, transitional, and quasi-periodic dynamics. Several single-parameter and multi-parameter training strategies were applied to train the reduced-order model. We found that \rev{an adequate} prediction of solutions beyond the set of parameters and initial conditions used to generate the training data requires including more snapshots that represent chaotic or short-time transient behaviors of the gKS model.

The results suggest that POD-ROM is capable of reproducing statistical properties of turbulent solutions, such as power spectra, and accurately recovers the solution over relatively short time intervals. This prediction time depends on the chosen training strategy. For the transitional regime, characterized by an initial period of (quasi-)chaotic behavior followed by a more laminar regime, a properly trained POD-ROM can recover persistent patterns of the laminar solution. These patterns emerge in the ROM solution with a phase shift. We conclude that for such systems, simply computing the norm of the error between FOM and ROM solutions may not be a sufficient assessment of ROM quality.

The quasi-periodic regimes proved to be the most resistant to effective POD-ROM simulations when initial conditions and parameter values were not included in the training sets. Interestingly, ROMs trained with more data from chaotic and transitional regimes yielded somewhat better results in these cases. Developing more universal projection-based ROMs that perform equally well across all regimes—from turbulent to laminar—may require employing parameter-specific ROM bases, potentially exploiting ideas of tensor-based ROMs. We plan to explore this research direction in future work.

%%%%%%%%%%%%%%%%%%%%%%%%%%%%%%%%%
\section*{Acknowledgment} The authors R.B.M. and M.O. were supported in part by the U.S. National Science Foundation under the award DMS-2309197.

\section*{Code}
Code for simulating the gKS equation and ROM with DEIM is located on GitHub \cite{binmizan2026gks}.

\appendix
%%%%%%%%%%%%%%%%%%%%%%%%%%%%%%%%%
\section{Initial Conditions}
Initial conditions in the simulations of the KS and the gKS equation are given by \eqref{eq:initial_condition}. \\
Initial conditions for Figure \ref{fig5} were generated with $J=8$ using 
$A_j=\{$3.0233, 3.5171, -7.3354, -9.6221, -1.9993, 5.3954, -6.7118, -4.5122$\}\times10^{-2}$ and
$\phi_j=\{$6.08, 3.0139, 6.2266, 3.1967, 4.983, 3.0123, 3.6808, 2.2941$\}$. \\
IC1 and IC2 in this paper were generated with $J=8$. Specific values of 
coefficients and phases for IC1 are 
$A_j^{(1)}=\{$0.1593, -8.4911, 5.6365, 5.4644, -1.8904, 6.2566, 1.7180, -2.1826$\}\times 10^{-2}$ and 
$\phi_j^{(1)} = \{$3.5963, 2.1939, 4.1857, 5.4722, 5.5467, 4.5229, 1.2151, 0.16661$\}$.
Values of these parameters for IC2 are 
$A_j^{(2)}=\{$-1.5013, 7.649, 1.891, 8.5255, 5.531, 2.4718, 6.5647, 9.6781$\}\times 10^{-2}$ and 
$\phi_j^{(2)} = \{$3.3396, 3.6174, 1.7775, 5.3945, 2.7183, 0.1965, 5.1264, 4.5449$\}$.

%%%%%%%%%%%%%%%%%%%%%%%%%%%%%%%%%%%%%%%%
%\bibliographystyle{siam}
%\bibliography{refs}

\end{document}